\documentclass[twoside,11pt,draft]{article}

\usepackage{amssymb}

\renewcommand{\labelenumi}{(\theenumi)}

\newtheorem{theorem}{Theorem}[section]
\newtheorem{proposition}[theorem]{Proposition}
\newtheorem{lemma}[theorem]{Lemma}
\newtheorem{definition}[theorem]{Definition}
\newtheorem{example}[theorem]{Example}
\newtheorem{examples}[theorem]{Examples}
\newtheorem{remark}[theorem]{Remark}

\newcommand{\btheorem}{\begin{theorem}}
\newcommand{\etheorem}{\end{theorem}}
\newcommand{\bprop}{\begin{proposition}}
\newcommand{\eprop}{\end{proposition}}
\newcommand{\blemma}{\begin{lemma}}
\newcommand{\elemma}{\end{lemma}}
\newcommand{\bdefi}{\begin{definition}\rm}
\newcommand{\edefi}{\ \hfill $\blacktriangle$ \end{definition}}
\newcommand{\defiendpage}{\end{definition}}
\newcommand{\nredefi}{\nopagebreak \begin{flushright} \vspace{-3mm}%
 $\blacktriangle$ \end{flushright} \end{definition}}
\newcommand{\bex}{\begin{example}\rm}
\newcommand{\eex}{\ \hfill $\blacktriangle$ \end{example}}
\newcommand{\bexs}{\begin{examples} \mbox{ } \rm \begin{enumerate}
\vspace{-3mm}}
\newcommand{\eexs}{\end{enumerate}  \hfill $\blacktriangle$  \end{examples}}
\newcommand{\bremark}{\begin{remark}\rm}
\newcommand{\eremark}{\ \hfill $\blacktriangle$ \end{remark}}
\newcommand{\bremarks}{\begin{remark}\rm   \mbox{ } \rm \begin{enumerate}
\vspace{-3mm}}
\newcommand{\eremarks}{\end{enumerate}  \hfill $\blacktriangle$  \end{remark}}

\newcommand{\qed}{\ \hfill \rule{2mm}{2mm}}
\newenvironment{proof}{\noindent \textsc{Proof~:}~}{\qed\medskip\par}
\newcommand{\bproof}{\begin{proof}}
\newcommand{\eproof}{\end{proof} \hfill \vspace{6mm}} 
\newcommand{\proofendpage}{\enlargethispage{1cm} \end{proof} } 

\newcommand{\bvgln}{$$\begin{array}{lcl}}
\newcommand{\evgln}{\end{array}$$}
\newcommand{\bvgl}{\begin{eqnarray*}}
\newcommand{\evgl}{\end{eqnarray*}}

\newcommand{\cst}{C\/$^{\ast}$}

\newcommand{\Z}{\mbox{$\mathbb Z$} }
\newcommand{\R}{\mbox{$\mathbb R$} }
\newcommand{\C}{\mbox{$\mathbb C$} }

\newcommand{\1}{{\rm 1 \hspace{-0.22em}\hbox{\rm l}}}

\newcommand{\eps}{\varepsilon}
\renewcommand{\to}{\rightarrow}
\newcommand{\cross}{\rtimes}

\newcommand{\hoed}{\widehat}
\newcommand{\VN}{\mbox{$\mathcal{M}$}}
\newcommand{\punt}{\, \cdot \,}
\newcommand{\bolleke}{\overline}
\newcommand{\dual}[2]{\langle #1 \mid #2 \rangle}

\renewcommand{\H}{\mbox{$\mathcal{H}$}}
\newcommand{\K}{\mbox{$\mathcal{K}$}}
\newcommand{\B}{\mbox{$\mathcal{B}$}}

\newcommand{\M}{\mbox{$\mathcal{M}$}}

\newcommand{\lijntjetussen}{\vspace{5mm}}
\newcommand{\lijntjeomhoog}{\vspace{-5mm}}

\begin{document}

\begin{center}
{\huge The multiplicative unitary as a basis for duality} \vspace{3ex}

\vspace{2ex}

{\large May 2002

\vspace{5ex}

Ann Maes

Economische Hogeschool Sint-Aloysius

Stormstraat 2

B-1000 Brussel, Belgium

ann.maes@prof.ehsal.be

\vspace{2ex}

Alfons Van Daele

Departement wiskunde

Katholieke Universiteit Leuven

Celestijnenlaan 200B

B-3001 Heverlee, Belgium

alfons.vandaele@wis.kuleuven.ac.be

}

\end{center}
\vspace{5ex}

\begin{abstract}
\noindent The classical duality theory associates to an abelian locally compact group a dual
companion. Passing to a non-abelian group, a dual object can still be defined, but it is no
longer a group.  The search for a broader category which should include both the groups and
their duals, points towards the concept of quantization.

\noindent Classically, the regular representation of a group contains the complete information
about the structure of this group and its dual.  In this article, we follow Baaj and Skandalis
and study duality starting from an abstract version of such a representation~: the
multiplicative unitary.  We suggest extra conditions which will replace the regularity and
irreducibility of the multiplicative unitary.  From the proposed structure of a "quantum group
frame", we obtain two objects in duality.  We equip these objects with certain group-like
properties, which make them into candidate quantum groups.

\noindent We consider the concrete example of the quantum $az+b$-group, and discuss how it
fits into this framework.  Finally, we construct the crossed product of a quantum group frame
with a locally compact group.

\end{abstract}

\section{Introduction}

The classical duality theory associates to an abelian locally compact group a dual group.  The
famous theorem of Pontryagin and Van Kampen then says that this dual contains the complete
information about the structure of the group~: its dual is again the original group. Several
authors (Tannaka, Krein, Tatsuuma) have generalized this result to non-abelian groups, showing
that to an arbitrary locally compact group one can still associate a dual object, from which
the original group can be recovered (see~\cite{hr}). But here the beautiful symmetry of the
Pontryagin duality is lost~: this dual object is no longer a group. And so the search started
for a broader category which would include both the groups and their duals.

It appeared that a quantization of the group concept would bring the solution of this problem.
The idea behind the process of ''quantization'' is the following. In a first step, instead of
thinking of a space as a set of elements, one considers some set of functions on this space.
For instance, let $G$ be a measure space, and consider the set $L^{\infty}(G)$ of all
essentially bounded complex functions on $G$. With pointwise multiplication and supremum norm
it becomes a von Neumann algebra. If $G$ is equipped with a (measurable) product, it can be
lifted to a $^*$-homomorphism $\Phi : L^{\infty}(G) \rightarrow L^{\infty}(G \times G)$ given
by
 $$(\Phi(f))(p,q)=f(pq)$$
for $f \in L^{\infty}(G)$ and $p,q \in G$.  In the process of quantization we now forget about
the underlying space $G$ and consider instead the pair $(L^{\infty}(G), \Phi)$ of a
commutative von Neumann algebra with a comultiplication. The next step then is to drop the
commutativity condition on the von Neumann algebra.

\bdefi
Let $M$ be a von Neumann algebra. We denote by $M \otimes M$ the von Neumann algebraic tensor
product of $M$ with itself. A \emph{comultiplication} on $M$ is a normal, injective
\mbox{$^*$-homomorphism} $\Phi : M \rightarrow M \otimes M$ such that $\Phi(\1) = \1 \otimes
\1$ and $\Phi$ satisfies the \emph{coassociativity property} $(\iota \otimes \Phi )\Phi =
(\Phi \otimes \iota)\Phi$.
\edefi

In this definition, and further in this work, we use $\iota$ to denote the identity map.

What extra conditions do we need on $(M,\Phi)$ to have a locally compact quantum group ? This
question has proved to be a lot more difficult. The first idea that comes to the mind is to
generalize the notion of unit and inverse. This brings a lot of problems, since the natural
candidates for such generalized unit and inverse are maps which are not well-behaved in
general.

\lijntjetussen

A first framework for quantization of locally compact groups, was the theory of Kac algebras
(see~\cite{es} for an overview).  It defines a category which contains the locally compact
groups and which allows a duality theory within the category.  Hence it provides a complete
answer to the duality problem formulated above.

However, at this point of time new questions have arisen, to which the Kac algebras cannot
give a full answer.  First of all, while studying examples of ''quantized groups'', sometimes
a structure was found which did not fit into the Kac algebra framework. It looked like the Kac
algebra structure was too restrictive. Moreover, the definition of a Kac algebra claims much
more structure then the definition of a locally compact group.  For instance, the existence of
a ''Haar measure'' is included in the axioms, whereas in the case of a classical group, this
is a very remarkable result of the theory.

For the particular case of compact quantum groups, Woronowicz~\cite{wor3} has developed a
theory which answers all these questions (see~\cite{notes} for an overview of this theory).
For the general case, the most mature theory at this moment is the one of Kustermans and
Vaes~\cite{kv2,kv3,kv4}. It realises a large part of the quantization project~: the category
defined is larger than the Kac algebra category, and its axioms are more simple.  But the
existence of the Haar measure is still included in the axioms, and at this moment a proof for
this existence still seems to be far away.

\lijntjetussen

We will follow a different strategy to study duality.  Let $G$ be a locally compact group.
Fix a left Haar measure on $G$ and let $\H=L^2(G)$ denote the Hilbert space of (classes of)
square integrable complex functions on $G$.  Then the Hilbert space tensor product $\H \otimes
\H$ is (isomorphic to) the Hilbert space $L^2(G \times G)$.  Let $W_G \in \B(\H \otimes \H)$
be the unitary operator defined by
 $$(W_G \xi)(p,q) = \xi(p,p^{-1}q)$$
for $\xi \in L^2(G \times G)$ and $p,q \in G$.  It was shown by Stinespring~\cite{stinespring}
for the unimodular, and by Takesaki~\cite{tak2} for the general case, that $W_G$ contains the
complete information about the algebraic and topological structure of the group $G$.  The
work~\cite{vh} of Vanheeswijck where she studied group duality starting from this operator
also is to be situated in this scheme. Kac~\cite{kac1,kac2} defined this operator in an
abstract way and proved that the structure of two unimodular Kac algebras (in his terminology,
''ring groups'') in duality can equivalently be described in terms of the associated unitary.
Baaj and Skandalis~\cite{bs} went further in the study of operators ''of this kind'' on an
abstract level.  To a so-called ''multiplicative unitary'' with certain extra conditions, they
associate two candidate quantum groups in duality.  In this work, we continue on this path,
and propose an alternative for the conditions of regularity and irreducibility which Baaj and
Skandalis introduced.

\lijntjetussen

The text is organized as follows.  We start in section 2 with a revision of the properties of
a multiplicative unitary, and illustrate these in the case of a classical group. These are not
new results, they are merely included for the convenience of the reader.  In the third
section, we introduce our framework and deduce the main results. At the end of this section,
we compare our setting with the one of Baaj and Skandalis~\cite{bs}. In section 4 we study the
example of the quantum $az+b$-group. Finally, in section 5, we construct the crossed product
with a locally compact group.

This article reflects the work done while the second author was preparing her Ph.D.\ thesis. A
more detailed text can be found in~\cite{thesis}.

\lijntjetussen

The story of locally compact quantum groups is told in the language of operator algebras.  For
a general introduction into the vocabulary and grammar of this language, we refer
to~\cite{dix1,kr,stratila,takesaki}.  We introduce some notations.

If $\H$ is a Hilbert space, we denote by $\B(\H)$ the von Neumann algebra of all bounded
linear operators on $\H$. As a general convention, we use the symbol $\odot$ for the algebraic
tensor product of vector spaces, and $\otimes$ for completed tensor products. If $\H_1$ and
$\H_2$ are Hilbert spaces, $\H_1 \otimes \H_2$ denotes the Hilbert space tensor product of
$\H_1$ and $\H_2$. We will denote by $\Sigma$ the flip operator $\Sigma: \H_1 \otimes \H_2
\rightarrow \H_2 \otimes \H_1$ given by $$\Sigma (\xi \otimes \eta)=\eta \otimes \xi$$ if $\xi
\in \H_1$ and $\eta \in \H_2$. The corresponding flip acting on operators will be denoted by
$\sigma$, so $$\sigma (x \otimes y) = \Sigma (x \otimes y) \Sigma = y \otimes x$$ for any $x
\in \B(\H_1)$ and $y \in \B(\H_2)$.

When $M_1$ and $M_2$ are von Neumann algebras, $M_1 \otimes M_2$ will denote the von Neumann
algebraic tensor product of $M_1$ and $M_2$.  The predual of a von Neumann algebra $M$ will be
denoted by $M_{\ast}$.

The identity operator on a Hilbert space will be denoted by $\1$.  We will frequently use the
\emph{leg numbering notation}~: If $W \in \B(\H \otimes \H)$, we define $W_{12}, W_{23},
W_{13} \in \B(\H \otimes \H \otimes \H)$ by $W_{12} = W \otimes \1$, $W_{23}=\1 \otimes W$ and
$W_{13} = (\1 \otimes \Sigma)(W \otimes \1)(\1 \otimes \Sigma)$. Accordingly, leg numbering
will be used for tensor products of more than three copies of $\H$, or for tensor products of
different Hilbert spaces.

For given vectors $\xi, \eta$ of a Hilbert space $\H$, we will denote by $\omega_{\xi, \eta}$
the linear functional on $\B(\H)$ such that $\omega_{\xi, \eta}(x)=\langle x \xi,\eta \rangle$
for any $x \in \B(\H)$.

If $G$ is a locally compact group, the set of all complex continuous functions on $G$ with
compact support will be denoted by $K(G)$. We will use the left Haar measure, and denote it by
$dq$.

\lijntjetussen
{\bf Acknowledgement}

Part of this work was done during the authors' stay at the NTNU in Trondheim. The authors
would like to thank Magnus Landstad and Chris Skau for their warm hospitality.

\section{The multiplicative unitary}
We begin by recollecting some material of Baaj and Skandalis~\cite{bs}.  This section does not
contain any new results, but is included for completeness and for convenience of the reader.

\bdefi
Let $\H$ be a Hilbert space. A unitary operator $W \in \B(\H \otimes \H)$ is said to be
\emph{multiplicative} if it satisfies the \emph{Pentagon equation}
 $W_{12} W_{13} W_{23} = W_{23}W_{12}.$
\edefi

\bex
Let $G$ be a locally compact group.  Let $W_G \in B(\H \otimes \H)$ be the operator defined by
 $$(W_G\, \xi) (p,q) = \xi (p,p^{-1} q)$$
for $\xi \in L^2(G \times G)$ and $p,q \in G $. Then $W_G$ is a multiplicative unitary~: the
Pentagon equation here amounts to the associativity of the product of $G$.
\eex

For the rest of this section, we fix a Hilbert space $\H$ and a multiplicative unitary $W$ on
$\H$. Define the sets
 $$ \begin{array}{ccccc}
 S &=& S(W) &=& \{( \iota \otimes \omega) (W) \mid \omega \in \B(\H)_{\ast} \},\\
 \hoed{S} &=& \hoed{S}(W) &=& \{(\omega \otimes \iota ) (W) \mid \omega \in \B(\H)_{\ast} \}.
 \end{array}$$

\begin{proposition} \textup{(\cite{bs} proposition 1.4)}
The vector spaces $S$ and $\hoed{S}$ are subalgebras of $\B(\H)$. They act non-degenerately on
$\H$.
\end{proposition}

Let $M = M(W)$ be the von Neumann algebra generated by $S$, i.e.\ the smallest von Neumann
algebra which contains $S$.  Let $\hoed{M} = \hoed{M}(W)$ be the von Neumann algebra generated
by $\hoed{S}$.

\begin{proposition} \label{prop-Win}
$W \in M \otimes \hoed{M}$.
\end{proposition}

\begin{proof}
Clearly $W$ commutes with $M' \otimes \1$ and with $\1 \otimes \hoed{M}'$, so $W \in (M'
\otimes \hoed{M}')'$.  By the commutation theorem for tensor products of von Neumann algebras
(see e.g.~\cite{takesaki} theorem\ 5.9), this last algebra equals $M \otimes \hoed{M}$.
\end{proof}

The comultiplication on $M$ will be provided by the multiplicative unitary. For $m \in M$ let
$\Phi(m) = W^{\ast}(\1 \otimes m)W$.  Then $\Phi$ maps $M$ into $M \otimes M$~: For $\omega
\in \B(\H)$ we have
 $$  W^{\ast}(\1 \otimes (\iota \otimes \omega )(W))W
 = (\iota \otimes \iota \otimes \omega )(W^{\ast}_{12} W_{23} W_{12})\\
 = (\iota \otimes \iota \otimes \omega)(W_{13} W_{23}).$$
Hence by proposition~\ref{prop-Win}, it follows that $\Phi(m) \in M \otimes M$ for any $m \in
S$. The set $\{m \in M \mid \Phi(m) \in M \otimes M \}$ is a von Neumann algebra. Since it
contains $S$, it must be the whole of $M$.

\begin{proposition}
The map $\Phi : M \rightarrow M \otimes M$ is a comultiplication.
\end{proposition}

\begin{proof}
The coassociativity is provided by the multiplicativity of $W$.  Indeed, for $m \in M$ we have
 \bvgl
 (\Phi \otimes \iota)\Phi (m)
 &=& W^{\ast}_{12}W^{\ast}_{23}(\1 \otimes \1 \otimes m) W_{23} W_{12}\\
 &=& W^{\ast}_{23}W^{\ast}_{13}W^{\ast}_{12}(\1 \otimes \1 \otimes m) W_{12}W_{13}W_{23}
 = (\iota \otimes \Phi)\Phi(m). \vspace{-0.7cm}
 \evgl
\end{proof}

\bremark
If $W \in \B(\H \otimes \H)$ is a multiplicative unitary, then also $\Sigma W^{\ast} \Sigma$
is. We will denote it by $\hoed{W}$ and call it the dual multiplicative unitary.  We write
$S(W)^{\ast}$ for the set $\{x^{\ast} \mid x \in S(W)\}$, and in the same way we use the
notation $\hoed{S}(W)^{\ast}$. Then we have that
  $$S(\hoed{W}) = \{(\iota \otimes \omega)
 (\Sigma W^{\ast} \Sigma) \mid \omega \in \B(\H)_{\ast}\} = \hoed{S}(W)^{\ast};$$
and similarly $\hoed{S}(\hoed{W}) = S(W)^{\ast}$.

The operator $W$ establishes a duality between $S$ and $\hoed{S}$~: For $\omega$, $\psi \in
\B(\H)_{\ast}$ let
 $$\langle (\iota \otimes \omega)(W) \mid
 (\psi \otimes \iota)(W) \rangle :=(\psi \otimes \omega)(W)
 = \omega ((\psi \otimes \iota)(W)).$$
This gives a well-defined pairing $\langle \punt \mid \punt \rangle$ on $S \times \hoed{S}$.
This pairing translates the comultiplication on $M$ into the product of $\hoed{M}$.  We give a
rigorous argument only in the case that $\H$ is a finite-dimensional Hilbert space. Let
$\{e_1, e_2, \ldots, e_n \}$ be an orthonormal basis for $\H$.  Let $\omega = \omega_{\xi,
\eta}$ with $\xi, \eta \in \H$, and let $\psi_1, \psi_2 \in \B(\H)_{\ast}$. With $m = (\iota
\otimes \omega)(W)$ and $\hoed{m_1} = (\psi_1 \otimes \iota)(W)$, $\hoed{m_2} = (\psi_2
\otimes \iota)(W)$, we have
 \bvgl
 \lefteqn{\langle \Phi(m) \mid \hoed{m_1} \otimes \hoed{m_2}\rangle =
 \langle (\iota \otimes \iota \otimes \omega) (W^{\ast}_{12}W_{23}W_{12}) \mid
 (\psi_1 \otimes \iota)(W) \otimes (\psi_2 \otimes \iota)(W) \rangle
 \mbox{\hspace{1cm}}} \hspace{1cm} \\
 &=& \langle (\iota \otimes \iota \otimes \omega) (W_{13}W_{23}) \mid
 (\psi_1 \otimes \iota)(W) \otimes (\psi_2 \otimes \iota)(W) \rangle
 \mbox{\hspace{1cm}} \\
 &=& \sum_{i=1}^{n} (\psi_1 \otimes \omega_{e_i, \eta}) (W)\,
 (\psi_2 \otimes \omega_{\xi, e_i})(W) \:
 = \: (\psi_1 \otimes \psi_2 \otimes \omega) (W_{13}W_{23})\\
 &=& \omega (\hoed{m_1} \hoed{m_2}) \:
 = \: \langle m \mid \hoed{m_1} \hoed{m_2}\rangle.
 \evgl
The dual multiplicative unitary $\hoed{W}$ gives a coproduct $\hoed{\Phi}$ on the dual von
Neumann algebra $\hoed{M}$, by $\hoed{\Phi}(\hoed{m}) = \hoed{W}^{\ast}(\1 \otimes \hoed{m})
\hoed{W}$.  Note that also the flipped map $\Phi' = \sigma \circ \Phi$ is a comultiplication
on $\hoed{M}$. With a similar calculation as the one above, one obtains that the product on
$M$ is dual to this {\it flipped coproduct} on $\hoed{M}$~:  For $m_1,m_2 \in M$ and $\hoed{m}
\in \hoed{M}$ we have $$\langle m_1 \otimes m_2 \mid \hoed{\Phi}'(\hoed{m})\rangle = \langle
m_1 m_2 \mid \hoed{m} \rangle.$$
\eremark

\bex \label{ex-group1}
Let $G$ be a locally compact group. Let $\lambda : G \to \B(L^2(G)):p \mapsto \lambda_p$
denote the \emph{left regular representation} of $G$, given by $(\lambda_p
\xi)(q)=\xi(p^{-1}q)$ if $\xi \in L^2(G)$ and $p,q \in G$.  The von Neumann algebra generated
by $\{ \lambda_{p} \mid p \in G \}$ is the \emph{group von Neumann algebra} $\M(G)$. It is a
well-known result (see~\cite{hr} theorem 22.11) that $\lambda$ also defines a representation
of $L^1(G)$ on $L^2(G)$.  We will denote this representation again by $\lambda$.  It is given
by
 $$(\lambda (f) \xi)(p) = \int f(q)(\lambda_q(\xi))(p)\, dq
 = \int f(q) \xi (q^{-1}p)\,dq$$
for $f \in K(G) \subseteq L^1(G)$, for $\xi \in K(G) \subseteq L^2(G)$ and $p \in G$.  The von
Neumann algebra generated by $\{\lambda(f) \mid f \in L^1(G)\}$ is again $\VN(G)$.

Let $W_G$ be the multiplicative unitary associated to $G$.  For $\omega =
\omega_{\eta_1,\eta_2}$ with $\eta_1, \eta_2 \in K(G)$, we have that
 $$((\iota \otimes \omega)(W_{G})\xi)(p)= \xi(p) \int \eta_1 (p^{-1}q)\:
 \overline{\eta_2(q)}\,dq $$
if  $\xi \in L^2(G)$ and $p \in G$. The von Neumann algebra generated by the functions
 $$G \rightarrow \C : p \mapsto \int \eta_1(p^{-1}q)\:
 \overline{\eta_2(q)}\,dq$$
(represented as multiplication operators on $L^2(G)$), is $L^{\infty} (G)$.  So we have that
$M(W_G) = L^{\infty}(G)$.  The comultiplication induced by $W_G$ is given by $\Phi(f)(p,q) =
f(pq)$ if $f \in L^{\infty}(G)$ and $p,q \in G$. Let $\eta_1, \eta_2 \in K(G)$ and denote the
function $\eta_1 \overline{\eta_2}$ in $K(G)$ by $f$. Then we have that
 \begin{equation} \label{eq-gvna}
 ((\omega_{\eta_1, \eta_2} \otimes \iota) (W_{G}) \xi)(p)= (\lambda (f) \xi)(p)
 \end{equation}
if $\xi \in L^2(G)$ and $p \in G$. It follows that $\hoed{M}(W_G) = \VN(G)$. The corresponding
comultiplication is given by $\hoed{\Phi}(\lambda_p) = \lambda_p \otimes \lambda_p$ for $p \in
G$.

So our group G gives rise to two ''quantum groups'' which are in duality~: On the one hand we
have the function algebra $L^{\infty} (G)$ where multiplication is trivial and
comultiplication reflects the group product. On the other hand there is the group von Neumann
algebra $\VN(G)$ where multiplication reflects the group product and comultiplication is
trivial.  We calculate the duality. Let $\omega \in \B(\H)_{\ast}$.  We saw that $(\iota
\otimes \omega)(W_G)$ is (the operator  of multiplication by) a function $f \in
L^{\infty}(G)$.  Let $\eta_1, \eta_2 \in K(G)$, and let $g = \eta_1 \overline{\eta_2} \in
K(G)$. Then equality (\ref{eq-gvna}) gives us that $\lambda(g) = (\omega_{\eta_1,\eta_2}
\otimes \iota)(W_G)$.  We have
 \bvgl
 \dual{f}{\lambda(g)}
 &=& \dual{(\iota \otimes \omega)(W_G)} {(\omega_{\eta_1,\eta_2} \otimes \iota)(W_G)}\\
 &=& (\omega_{\eta_1,\eta_2} \otimes \omega)(W_G) = \omega_{\eta_1,\eta_2}(f)\\
 &=& \int f(p)g(p) \, dp.
 \evgl

Let us have a look at the case where $G$ is an abelian group.  Denote the Pontryagin dual of
$G$ by $\hoed{G}$, and the duality between $G$ and $\hoed{G}$ by $\dual{\punt}{\punt}$.  By
Plancherel's theorem, for the suitable choice of the Haar measures on $G$ and $\hoed{G}$, the
Fourier transform $L^2(G) \to L^2(G):\xi \to \hoed{\xi}$ is an isomorphism of Hilbert spaces.
It translates the dual multiplicative unitary $\hoed{W_G}$ to the multiplicative unitary
$W_{\hoed{G}}$ corresponding to the dual group.  Indeed, for $\xi \in K(G \times G)$ and
$\gamma, \nu \in \hoed{G}$ we have
 \bvgl
 (\hoed{W_G}\xi)\hoed{\hspace{0.5em}}(\gamma,\nu)
 &=& \int \int \dual{\gamma}{p}^- \, \dual{\nu}{q}^- \, \xi(qp,q) \,dp \,dq \\
 &=& \int \int \dual{\gamma}{q^{-1}p}^- \, \dual{\nu}{q}^- \, \xi(p,q) \,dp \,dq \\
 &=& \int \int \dual{\gamma}{p}^- \, \dual{\gamma^{-1}}{q}^- \,
 \dual{\nu}{q}^- \, \xi(p,q) \,dp \,dq \\
 &=& \hoed{\xi}(\gamma, \gamma^{-1}\nu) = (W_{\hoed{G}} \hat{\xi})(\gamma,\nu).
 \evgl
In this way the theory of multiplicative unitaries generalizes the classical Pontryagin
duality.

We would like to remark that this is a very well-known example.  Nevertheless we have chosen
to include it because of it's didactic value.
\eex

\section{Quantum Group Frames}

When the von Neumann algebra associated to a multiplicative unitary is commutative, this
unitary arises from a classical locally compact group, which can be recovered completely from
it (\cite{bs}, section 2). Unfortunately, the general quantum case is not so straightforward.
In order to deduce a number of "group-like" properties, in this section we introduce some
extra structure.  At the end of the section, we compare this structure with the conditions of
Baaj and Skandalis in~\cite{bs}.

\lijntjetussen
Note that if $\H$ and $\K$ are Hilbert spaces, and $W \in \B (\H \otimes \H)$ is a
multiplicative unitary, then $W_{13}=(\iota \otimes \sigma \otimes \iota) (W \otimes \1_{\K}
\otimes \1_{\K})$ is a multiplicative unitary on $\H \otimes \K \otimes \H \otimes \K$. Of
course this operator will not really carry more information than $W$ itself. To avoid this
situation, we introduce the notion ''trim''.  Let $W$ be a multiplicative unitary, and let
$S$, $S^{\ast}$, $\hoed{S}$, $\hoed{S}^{\ast}$, $M$ and $\hoed{M}$ be as before. 
We denote by $S \hoed{S}^{\ast}$ the linear span of the set $\{s \hoed{s}^{\ast} \mid s \in S
\mbox{ and } \hoed{s} \in \hoed{S} \}$.

\begin{proposition} \label{prop-trim}
The $\sigma$-weak closure of $S \hoed{S}^{\ast}$ is a subalgebra of $\B(\H)$.
\end{proposition}

\begin{proof}
We denote the $\sigma$-weak closure of a set $X$ by $X^-$. For $\omega, \omega' \in
\B(\H)_{\ast}$ we have
 \bvgl
 (\omega \otimes \iota)(W^{\ast}) \: (\iota \otimes \omega')(W)
 &=& (\omega \otimes \iota \otimes \omega')(W_{12}^{\ast}W_{23})\\
 &=& (\omega \otimes \iota \otimes \omega')(W_{13} W_{23} W_{12}^{\ast}).
 \evgl

As $W \in M \otimes \hoed{M}$, it can be $\sigma$-weakly approximated by elements of the form
$\sum_{i=1}^{n}m_i \otimes \hoed{m_i}$ with $m_i \in M$ and $\hoed{m_i} \in \hoed{M}$. Using
this we obtain that the operator above can be $\sigma$-weakly approximated by elements
 $$ \sum_{i=1}^{n} (\iota \otimes \omega' (\hoed{m_i}\punt))(W) \:
 (\omega(m_{i}\punt ) \otimes \iota)(W^{\ast}) $$
of $S \hoed{S}^{\ast}$. So $ \hoed{S}^{\ast}S$ is contained in $(S \hoed{S}^{\ast})^-$. Hence
$(S \hoed{S}^{\ast})^- (S \hoed{S}^{\ast})^-$ is contained in $(S S \hoed{S}^{\ast}
\hoed{S}^{\ast})^- \subseteq (S \hoed{S}^{\ast})^-$.
\end{proof}

\bdefi \label{deftrim}
We call the multiplicative unitary $W$ \emph{trim} if the algebra $S \hoed{S}^{\ast}$ is
$\sigma$-weakly dense in $\B(\H)$. \edefi

We have that the algebra $S \hoed{S}^{\ast}$ acts non-degenerately on $\H$. Hence {\it if the
\mbox{$\sigma$-weak} closures of $S$ and $\hoed{S}$ are $^*$-algebras}, the $\sigma$-weak
closure of $S \hoed{S}$ will be a von Neumann algebra.  When this von Neumann algebra is not
the whole of $\B(\H)$, there is a non-trivial projection $p \in \B (\H)$ which commutes with
both $S$ and $\hoed{S}$. In that case we can "cut down $W$" by this projection, and study the
multiplicative unitary $W(p \otimes p)=(p \otimes p)W(p \otimes p)$ on the restricted Hilbert
space $p \H \otimes p\H$.  It is clear that this restricted multiplicative unitary still
carries essentially all the information of the underlying duality structure.

\bex \label{ex-group2}
Let $G$ be a locally compact group, and $W_G$ be as in example~\ref{ex-group1}. Since
$L^{\infty}(G)' = L^{\infty}(G)$, we have that $L^{\infty}(G)' \cap \VN(G)' = \C \1$. Hence
$L^{\infty}(G) \M(G)$ is $\sigma$-weakly dense in $\B(\H)$ and $W_G$ is trim.
\end{example}

Now we are ready to introduce our framework.

\bdefi \label{def-qgf}
Let $\H$ be a Hilbert space. A \emph{quantum group frame} on $\H$ is a triple $(W,J,\hoed{J})$
such that
\begin{enumerate}
\item $W \in \B ( \H \otimes \H )$ is a
multiplicative unitary;
\item $W$ is trim (cf.\ definition~\ref{deftrim});
\item $J$ and $\hoed{J}$ are anti-linear
operators on $\H$ such that $J= J^{\ast}$, $J^{2} = \1$, $\hoed{J} = \hoed{J}^{\ast}$,
$\hoed{J}^{2}= \1$ and $W^{\ast} = ( \hoed{J} \otimes J)W( \hoed{J} \otimes J)$; \nopagebreak
\item $J M J \subseteq M'$ and $\hoed{J} \hoed{M} \hoed{J}
\subseteq \hoed{M}'$.
\end{enumerate} \lijntjeomhoog
\edefi

\bremark
The definition of a quantum group frame is self-dual. Indeed, $S(\hoed{W})
\hoed{S}(\hoed{W})^{\ast}$ is the algebra $\hoed{S}(W)^{\ast}S(W)$; by the argument in the
proof of proposition~\ref{prop-trim} it has the same closure as the algebra
$S(W)\hoed{S}(W)^{\ast}$. Therefore $\hoed{W}$ is trim if and only if $W$ is trim. And if
$(W,J,\hoed{J})$ is a quantum group frame then also $(\hoed{W}, \hoed{J}, J)$ is one.
\eremark

We postpone the discussion about the meaning of the operators $J$ and $\hoed{J}$ for a little
while, and first show that the $\sigma$-weak closures of the algebras $S$ and $\hoed{S}$
associated to a quantum group frame are $^*$-algebras --- see theorem~\ref{th-selfadjoint}
below. Hence, since $S$ and $\hoed{S}$ act non-degenerately on $\H$, these closures will
coincide with the von Neumann algebras $M$ and $\hoed{M}$.  This property is crucial for the
further development of the theory.  Although we don't know of any examples where it is not
satisfied, it is not a straightforward result.  In order to obtain the corresponding result on
the \cst-algebraic level (i.e. to deduce that the norm-closure of $S$ is a $^*$-algebra), Baaj
and Skandalis \cite{bs} introduced the notion of regularity.  This notion later was weakened
by Baaj in his work~\cite{baaj}.  We will compare our axioms with their notions in
proposition~\ref{prop-equivalent}.

\lijntjetussen
We start with some technical results.

\begin{lemma} \label{lemma-predual}
Let $\H$ be a Hilbert space, and let $B$ be a linear subspace of $\B(\H)$ such that the weak
closure of $B$ contains $\1$.  Then the vector spaces generated by
 \bvgl
 &&\{\omega (\punt b) \mid \omega \in \B(\H)_{\ast},\ b \in B \}\\
 &&\{\omega (b \punt ) \mid \omega \in \B(\H)_{\ast},\ b \in B \},
 \evgl
are both norm dense subspaces of $\B(\H)_{\ast}$.
\end{lemma}

\begin{proof}
We will prove the first vector space to be norm dense in $\B(\H)_{\ast}$.  The norm density of
the second one follows analogously.

Since the weak closure of $B$ contains $\1$, we have that $B \H$ is weakly dense in $\H$.
Hence, because it is a subspace, $B \H$ is also norm dense in $\H$.

From this, it follows that the closed vector space generated by $$\{\omega_{b \xi,\eta} \mid b
\in B;\ \xi,\eta \in \H \}$$ in $\B(\H)^{\ast}$ is the same as the one generated by
$\{\omega_{\xi,\eta} \mid  \xi,\eta \in \H \}.$  The latter subset generates $\B(\H)_{\ast}$.
\end{proof}

\lijntjetussen

\begin{lemma} \label{lemma-st}
Let $W$ be a multiplicative unitary on a Hilbert space $H$, and $J, \hoed{J}$ anti-linear
operators on $\H$ such that $J= J^{\ast}$, $J^{2} = \1$, $\hoed{J} = \hoed{J}^{\ast}$,
$\hoed{J}^{2}= \1$ and $W^{\ast} = ( \hoed{J} \otimes J)W( \hoed{J} \otimes J)$.  Then
$\hoed{J}S\hoed{J}=S^{\ast}$ and $J \hoed{S} J = \hoed{S}^{\ast}$.
\end{lemma}

\begin{proof}
For $\omega \in \B(\H)_{\ast}$ we have that
 $$\hoed{J}(\iota \otimes \omega)(W)\hoed{J}
 = (\iota \otimes \theta)(W^{\ast}),$$
where $\theta$ is the linear functional on $\B(\H)$ given by $\theta(x)=\omega(JxJ)^-$ for any
$x \in \B(\H)$.  Hence $\hoed{J}S\hoed{J}=S^{\ast}$. The second result follows analogously.
\end{proof}

 \lijntjetussen

Now let $(W,J,\hoed{J})$ be a quantum group frame on a Hilbert space $\H$.

\begin{lemma} \label{lemma-set}
The linear span of the set $$\{(x \otimes \1)W(\1 \otimes y) \mid x,y \in \B(\H)\}$$ is
$\sigma$-weakly dense in $\B(\H) \otimes \B(\H)$.
\end{lemma}

\begin{proof}
Let $\xi, \eta \in \H$ and let $\omega$ denote the vector functional $\omega_{\xi, \eta}$. The
multiplicativity of $W$ gives us
 $$ \1 \otimes (\iota \otimes \omega)(W)
 = (\iota \otimes \iota \otimes
 \omega)(W_{13}^{\ast}W_{12}^{\ast}W_{23}W_{12}).$$
Let $\{e_i \mid i \in I \}$ be an orthonormal basis for $\H$.  Let $p_i$ denote the projection
operator on $\C e_i$. Then $\sum_{i \in I} p_i = \1$, where we have convergence in (e.g.)\ the
$\sigma$-weak topology.  Now
 \bvgl
 \1 \otimes (\iota \otimes \omega)(W)
 &=& \sum_{i \in I} (\iota \otimes \iota \otimes \omega)(W^{\ast}_{13}
 (\1 \otimes \1 \otimes p_i) W_{12}^{\ast}W_{23}W_{12})\\
 &=& \sum_{i \in I}
 (\iota \otimes \iota \otimes \omega_{e_i,\eta})(W^{\ast}_{13}) \: W^{\ast} \:
 (\iota \otimes \iota \otimes \omega_{\xi,e_i})(W_{23}) \: W,
 \evgl
with convergence in the $\sigma$-weak topology.   We will denote the $\sigma$-weak closure of
the linear span of a set $X$ by $\overline{\mbox{sp}}(X)$.  The above calculation gives us
that
 $$\1 \otimes S \subseteq
 \overline{\mbox{sp}}\{(s^{\ast} \otimes \1)W^{\ast}(\1 \otimes s')W \mid s,s' \in S \}.$$
Multiplying these sets by the vector space $\B(\H) \odot \hoed{M}'$ on the left hand side
(recall that $\odot$ denotes the algebraic tensor product of vector spaces) gives
 $$\B(\H) \odot \hoed{M}'S \subseteq \overline{\mbox{sp}}
 \{(x \otimes \1)W^{\ast}(\1 \otimes y)W \mid x,y \in \B(\H)\}.$$
Axiom (2) of definition~\ref{def-qgf} gives us that the set $S\hoed{S}^{\ast}$ is
$\sigma$-weakly dense in $\B(\H)$, so also $\hoed{J}S\hoed{J}\hoed{J}\hoed{S}^{\ast}\hoed{J}$
is $\sigma$-weakly dense in $\B(\H)$. By lemma~\ref{lemma-st} and axiom (4), the latter set is
contained in $S^{\ast}\hoed{M}'$, which therefore is dense in $\B(\H)$ as well. The same then
holds for its adjoint $\hoed{M}'S$.  So we obtain that
 $$\B(\H) \otimes \B(\H)
 = \overline{\mbox{sp}} \{(x \otimes \1)W^{\ast}(\1 \otimes y)W
 \mid x,y \in \B(\H)\}.$$
Multiplying these spaces on the right hand side with the unitary $W^{\ast}$, we have
 $$\B(\H) \otimes \B(\H) = \overline{\mbox{sp}} \{(x \otimes \1)
 W^{\ast}(\1 \otimes y) \mid x,y \in \B(\H)\}.$$
Using $W^{\ast}=(\hoed{J} \otimes J)W(\hoed{J} \otimes J)$ now gives the result.
\end{proof}

\lijntjetussen
Following Baaj and Skandalis \cite{bs}, we introduce the vector space
 \begin{equation} \label{eq-C}
 C = \{ (\iota \otimes \omega)(\Sigma W)\mid \omega \in \B ( \H)_{\ast} \}.
 \end{equation}
Similarly, let $D$ be the vector space
 \begin{equation}
 D = \{(\omega \otimes \iota)(\Sigma W) \mid \omega \in \B(\H)_{\ast} \}
 = \{(\iota \otimes \omega)(W \Sigma) \mid \omega \in \B(\H)_{\ast} \}.
 \label{eq-D}
 \end{equation}

\begin{lemma} \label{lemma-C}
The vector spaces $C$ and $D$ are $\sigma$-weakly dense in $\B(\H)$.
\end{lemma}

\begin{proof}
We give the proof for $C$.  The result for $D$ is proved analogously, or follows from the
previous by using that $C^{\ast}=JD\hoed{J}$.  Let $\xi,\xi',\eta,\eta' \in \H$ and let $x,y$
be the rank one operators $x=\langle \punt,\xi \rangle \xi'$ and $y= \langle \punt ,\eta
\rangle \eta'$. We have that
 $$ \Sigma (x \otimes \1) W (\1 \otimes y)
 = (\1 \otimes x) \Sigma W (\1 \otimes y)\\
 = (\iota \otimes \omega_{\eta',\xi})(\Sigma W) \otimes
 \langle \punt , \eta \rangle \xi'.$$
By lemma~\ref{lemma-set}, the closure of the vector space spanned by
 $$\{(\iota \otimes \omega) (\Sigma (x \otimes \1) W (\1 \otimes y))
 \mid \omega \in \B(\H)_{\ast}; \: x,y \mbox{ are rank one operators on }\H\}$$
is $\B(\H)$. Hence we obtain that the closure of $C$ is $\B(\H)$.
\end{proof}

We then can use the technique of Baaj and Skandalis (\cite{bs} proposition 3.5) to obtain the
following result.

\begin{theorem} \label{th-selfadjoint}
Let $(W,J,\hoed{J})$ be a quantum group frame. Then the $\sigma$-weak closure of $S$ is a
$^*$-algebra. Hence it coincides with the von Neumann algebra $M$. Similarly, the
$\sigma$-weak closure of $\hoed{S}$ is the von Neumann algebra $\hoed{M}$.
\end{theorem}

\begin{proof}
As before, we denote the $\sigma$-weak closure of the linear span of a set $X$ by
$\overline{\mbox{sp}}(X)$. Since by lemma~\ref{lemma-C}, the set $D$ is $\sigma$-weakly dense
in $\B(\H)$, it follows from lemma~\ref{lemma-predual} that the set $\{\omega'(z \punt) \mid z
\in D, \omega \in \B(\H)_{\ast}\}$ is dense in $\B(\H)_{\ast}$. We have
 \bvgl
 (S^{\ast})^- &=& \overline{\mbox{sp}}\{(\iota \otimes \omega( \punt z)) (W^{\ast})
 \mid  z=(\iota \otimes \omega')(W \Sigma); \omega, \omega' \in \B(\H)_{\ast} \}\\
 &=& \overline{\mbox{sp}}\{(\iota \otimes \omega \otimes \omega')
 (W^{\ast}_{12}W_{23} \Sigma_{23}) \mid \omega, \omega' \in \B(\H)_{\ast} \}\\
 &=& \overline{\mbox{sp}}\{(\iota \otimes \omega \otimes \omega')
 (W_{13}W_{23}W^{\ast}_{12} \Sigma_{23}) \mid \omega, \omega' \in  \B(\H)_{\ast} \}\\
 &=& \overline{\mbox{sp}}\{(\iota \otimes \omega') (W (\1 \otimes y) W^{\ast})
 \mid y = (\omega \otimes \iota)(W \Sigma);\ \omega, \omega' \in \B(\H)_{\ast} \}\\
 &=& \overline{\mbox{sp}}\{(\iota \otimes \omega') (W (\1 \otimes y) W^{\ast})
 \mid \omega' \in \B(\H)_{\ast}, y \in C \}.
 \evgl
Since the $\sigma$-weak closure of $C$ is $\B(\H)$, it follows that $(S^{\ast})^-$ is
self-adjoint, so also $S^-$ is self-adjoint.  By duality we obtain that also the $\sigma$-weak
closure of $\hoed{S}$ is self-adjoint.
\end{proof}

 \lijntjetussen

Let us now have a closer look at the operators $J$ and $\hoed{J}$.  We first return to the
example~\ref{ex-group2} of a classical group $G$. Let $J_G$  and $\hoed{J_G}$ be the
anti-linear operators on $L^{2} (G)$ defined by
 \bvgl
 (J_G \xi) (p) &=& \overline{\xi (p)}\\
 (\hoed{J_G} \xi) (p) &=& \Delta (p)^{-\frac{1}{2}} \, \overline{\xi (p^{-1}}),
 \evgl
where $\Delta$ denotes the modular function of $G$. It is straightforward to check that axioms
(3) and (4) are satisfied. We argued in example~\ref{ex-group2} that also axiom (2) is
satisfied.  Hence $(W_G, J_G, \hoed{J_G} )$ is a quantum group frame.

 \lijntjetussen

So in the group example, we see that $\hoed{J}$ is linked with the inverse operation. In the
quantum world, the role of the inverse is played by the antipode. In general this antipode
will not be a bounded operator, but $\hoed{J}$ will implement its unitary part (and so will
$J$ for the dual picture) --- see~\cite{kv2}. Note that by lemma~\ref{lemma-st}, we have that
$\hoed{J} M \hoed{J} = M$ and $J \hoed{M} J= \hoed{M}$. So there are linear maps
 $$ \begin{array}{ccccc}
 R &:& M \rightarrow M &:& m \mapsto \hoed{J} m^{\ast} \hoed{J},\\
 \hoed{R} &:& \hoed{M} \rightarrow \hoed{M} &:& \hoed{m} \mapsto J \hoed{m}^{\ast} J.
 \end{array} $$
In the group case, $R: L^{\infty}(G) \rightarrow L^{\infty}(G)$ is given by
$(Rf)(p)=f(p^{-1})$ whenever $f \in L^{\infty}(G)$, and $\hoed{R}: \M(G) \rightarrow \M(G)$ is
given by $\hoed{R}(\lambda_p)=\lambda_{p^{-1}}$ for $p \in G$.

 \lijntjetussen

The maps $R$ and $\hoed{R}$ are $^*$-anti-automorphisms.  Moreover, they flip the
comultiplication~:

\begin{proposition} \label{Rflips}
We have that
 \bvgl
 \Phi \circ R &=& \sigma \circ (R \otimes R) \circ \Phi, \\
 \hoed{\Phi} \circ \hoed{R} &=& \sigma \circ (\hoed{R} \otimes \hoed{R} ) \circ \hoed{\Phi}.
 \evgl
\end{proposition}

\begin{proof}
We have
 \bvgl
 W_{12}^{\ast} W_{23} W_{12}
 &=& W_{13} W_{23}= ( \hoed{J} \otimes \hoed{J} \otimes J)
 W_{13}^{\ast} W_{23}^{\ast} ( \hoed{J} \otimes \hoed{J} \otimes J)\\
 &=& ( \sigma \otimes \iota )(( \hoed{J} \otimes \hoed{J} \otimes J)
 W_{23}^{\ast} W_{13}^{\ast} ( \hoed{J} \otimes \hoed{J} \otimes J))\\
 &=& (\sigma \otimes \iota )(( \hoed{J} \otimes \hoed{J} \otimes J) W_{12}^{\ast}
 W_{23}^{\ast} W_{12} ( \hoed{J} \otimes \hoed{J} \otimes J)).
 \evgl
Let $\omega \in \B (\H)_{\ast}$ and let $m=( \iota \otimes \omega )(W)$. Define $\theta \in
\B(\H)_{\ast}$ by $\theta(x)=\omega (J x^{\ast} J)$. So $(\iota \otimes \theta)(W) = R(m)$.
Applying $( \iota \otimes \iota \otimes \theta)$ to the above equality gives $\Phi (R(m))=
\sigma ((R \otimes R) \Phi (m))$.  By theorem~\ref{th-selfadjoint}, such elements $m$ form a
$\sigma$-weakly dense part of the von Neumann algebra $M$. So the first statement follows. The
other equality follows by duality.
\end{proof}

 \lijntjetussen

In the following, we will discuss how our extra conditions on $W$ relate to the conditions of
regularity and irreducibility which Baaj and Skandalis impose in \cite{bs}. Let us fix a
Hilbert space $\H$ and a multiplicative unitary $W$ on $\H$.

As before, let $C=\{(\iota \otimes \omega)(\Sigma W) \mid \omega \in \B(\H)_{\ast}\}$.
Without extra conditions lemma~\ref{lemma-C} is not valid anymore.  But we have the following
property~:

\begin{lemma} \textup{(\cite{bs} proposition 3.2)} \newline
$C$ is an algebra which acts non-degenerately on $\H$.
\end{lemma}

Baaj and Skandalis~\cite{bs} define $W$ to be \emph{regular} if the norm closure of $C$ equals
$\K(\H )$, the \mbox{\cst-algebra} of compact operators on $\H$. Baaj \cite{baaj} defines $W$
to be \emph{semi-regular} if the norm closure of $C$ contains $\K ( \H )$. The other extra
condition which is imposed on $W$ in \cite{bs} is the irreducibility. The multiplicative
unitary $W$ is called \emph{irreducible} if there exists a unitary $U \in \B(\H )$ such that
\begin{enumerate}
\item $U^{2}=1$ and $(\Sigma ( \1 \otimes U)W)^{3}=1$,
\item the unitaries $\bolleke{W}= \Sigma (U \otimes \1 )W(U \otimes \1 )\Sigma$
and $\widetilde{W} = (U \otimes U) \bolleke{W} ( U \otimes U )$ are multiplicative.
\end{enumerate}

In our setting, the product $\hoed{J} J$ will play the role of $U$. Putting $U= \hoed{J} J$,
it is not immediately clear whether $U^{2}=1$, so we will define $\bolleke{W}$ and
$\widetilde{W}$ as
 $$ \begin{array}{ccccc}
 \bolleke{W} &=& \Sigma ( U^{\ast} \otimes \1 )W(U \otimes \1 ) \Sigma
 &=& (J \otimes J) \Sigma W^{\ast} \Sigma (J \otimes J)\\
 \widetilde{W} &=&(U \otimes U) \bolleke{W} ( U^{\ast} \otimes U^{\ast})
 &=& ( \hoed{J} \otimes \hoed{J} ) \Sigma W^{\ast} \Sigma ( \hoed{J} \otimes \hoed{J} ).
 \end{array} $$

\begin{proposition} \label{prop-equivalent}
Let $\H$ be a Hilbert space and $W \in \B ( \H \otimes \H )$ a multiplicative unitary. Suppose
that $J$ and $\hoed{J}$ are anti-linear operators on $\H$ and that $J= J^{\ast}$, $J^{2} =
\1$, $\hoed{J} = \hoed{J}^{\ast}$, $\hoed{J}^{2} = \1$ and $W^{\ast} = ( \hoed{J} \otimes J)W(
\hoed{J} \otimes J)$. Then the following two sets of conditions are equivalent~:
\begin{enumerate}
\item $W$ is trim,
\item $JMJ \subseteq M'$ and $\hoed{J} \hoed{M} \hoed{J} \subseteq
\hoed{M}'$;
\end{enumerate}
and
\renewcommand{\labelenumi}{(\theenumi')}
\begin{enumerate}
\item   $C$ is $\sigma$-weakly dense in $\B ( \H )$,
\item $(\1 \otimes \hoed{J} J) \Sigma \bolleke{W} W \widetilde{W}$ is
a scalar multiple of $\1$.
\end{enumerate}
\renewcommand{\labelenumi}{(\theenumi)}
\end{proposition}

Condition (1') is weaker than regularity and semi-regularity.   Of course, the existence of
the two anti-linear operators $J$ and $\hoed{J}$ is a stronger claim than the existence of
only their product. But from the examples it will become clear that these operators play a
very natural role, so that it is reasonable to work with $J$ and $\hoed{J}$ instead of their
product.  In this setting, the unitaries $\bolleke{W}$ and $\widetilde{W}$ are automatically
multiplicative.  This is another reason why these are natural conditions to impose.  Condition
(2') slightly weakens the further irreducibility condition $(\Sigma ( \1 \otimes U)W)^{3}=\1$,
which is equivalent to $(\1 \otimes U) \Sigma \bolleke{W} W \widetilde{W} = \1$.  In the
example of the quantum $az+b$-group (see section~\ref{sec-az+b}) we will indeed obtain a
scalar which is not $1$. Nevertheless, this turns out not to be really essential, as we will
see in the discussion after proposition~\ref{commute}.

Altogether, we think our conditions are more natural than conditions like (1') and (2'), and
therefore easier to verify in examples. However, a drawback is that we have not been able to
deduce that they imply the norm closure of $S$ to be a \cst-algebra.  But we did obtain that
its $\sigma$-weak closure is a von Neumann algebra.

\lijntjetussen
\begin{proof}
Suppose that the conditions (1) and (2) are satisfied.
\begin{itemize}
\item
Lemma \ref{lemma-C} gives us that condition (1') then holds.

\item We will show that $(\1 \otimes \hoed{J} J)\Sigma \bolleke{W} W
\widetilde{W}$ commutes with $M \otimes \1$. Similar calculations show that it also commutes
with $\hoed{M} \otimes \1$, $\1 \otimes M$ and $\1 \otimes \hoed{M}$. Hence by (1), it must be
a scalar multiple of $\1$.

First note that using the multiplicativity of $W$ we have
 \bvgl
 (\hoed{J} \otimes \hoed{J} \otimes J)
 W_{12}^{\ast} W_{23}^{\ast} (\hoed{J} \otimes \hoed{J} \otimes J)
 &=& (\hoed{J} \otimes \hoed{J} \otimes J)
 W_{23}^{\ast} W_{13}^{\ast} W_{12}^{\ast}
 (\hoed{J} \otimes \hoed{J} \otimes J)\\
 &=& \Sigma_{12} ( \hoed{J} \otimes \hoed{J} \otimes J)
 W_{13}^{\ast} W_{23}^{\ast} \Sigma_{12} W_{12}^{\ast}
 ( \hoed{J} \otimes \hoed{J} \otimes J).
 \evgl
Replacing $W^{\ast}$ by $(\hoed{J} \otimes J)W( \hoed{J} \otimes J)$ we obtain that
 $$
 (\1 \otimes \hoed{J} J \otimes \1) W_{12}
 (\1 \otimes J \hoed{J} \otimes \1 ) W_{23}
 = \Sigma_{12} W_{13} W_{23} ( \hoed{J} J \otimes \1 \otimes \1) \Sigma_{12}
 W_{12} ( \1 \otimes J \hoed{J} \otimes \1)$$
and using the multiplicativity once again, this gives
 $$
 (\1 \otimes \hoed{J} J \otimes \1) W_{12}
 (\1 \otimes J \hoed{J} \otimes \1 ) W_{23}
 = \Sigma_{12} W_{12}^{\ast} W_{23} W_{12}
 ( \hoed{J} J \otimes \1 \otimes \1) \Sigma_{12} W_{12}
 ( \1 \otimes J \hoed{J} \otimes \1).
 $$
Therefore, moving some of the factors of this equality of operators to the other side, we have
 \begin{eqnarray*}
 \lefteqn{W_{23}^{\ast} W_{12} ( \hoed{J} J \otimes \1 \otimes \1)
 \Sigma_{12} W_{12} ( \1 \otimes J \hoed{J} \otimes \1)} \hspace{4cm} \\
 & & = W_{12} (\hoed{J} J \otimes \1 \otimes \1 ) \Sigma_{12}
 W_{12} ( \1 \otimes J \hoed{J} \otimes \1) W_{23}^{\ast}.
 \end{eqnarray*}
It follows that $( \Sigma W( \hoed{J} J \otimes \1)\Sigma W( \1 \otimes J
\hoed{J})\Sigma)_{12}$ commutes with $W_{13}$, and hence the algebra $M \otimes \1$ commutes
with $\Sigma W( \hoed{J} J \otimes \1 ) \Sigma W( \1 \otimes J \hoed{J} ) \Sigma$. Using the
fact that $\hoed{J} M \hoed{J} = M$ and assumption (2) $JMJ \subseteq M'$, we may conclude
that $M \otimes \1$ commutes with $$(J \hoed{J} \otimes \hoed{J} J)W(\hoed{J} J \otimes
\1)\Sigma W( \hoed{J} J \otimes \1) \Sigma W(\1 \otimes J \hoed{J} )\Sigma = (\1 \otimes
\hoed{J} J) \Sigma \bolleke{W} W \widetilde{W}.$$
\end{itemize}

Conversely, suppose that $( \1 \otimes \hoed{J} J) \Sigma \bolleke{W} W \widetilde{W} = k \1
\otimes \1$ for some $k \in \C$.
\begin{itemize}
\item
Following the proof of proposition \ref{Rflips} we have
 \bvgl
 W_{12}^{\ast} W_{23} W_{12}
 &=& \Sigma_{12} ( \hoed{J} \otimes \hoed{J} \otimes J)
 W_{12}^{\ast} W_{23}^{\ast} W_{12} ( \hoed{J} \otimes \hoed{J} \otimes J) \Sigma_{12} \\
 &=& \Sigma_{12} ( \hoed{J} \otimes \hoed{J} \otimes J)
 W_{12}^{\ast} (J \otimes \hoed{J} \otimes J) W_{23} (J \otimes
 \hoed{J} \otimes {J}) W_{12} ( \hoed{J} \otimes \hoed{J} \otimes J) \Sigma_{12}.
 \evgl
Commuting $\Sigma_{12}$ with $(\hoed{J} \otimes \hoed{J} \otimes J)$ and moving the factors
$W_{12}$ of the left hand side and $( \hoed{J} \otimes \hoed{J} \otimes J) \Sigma_{12}
 W_{12}^{\ast} (J \otimes \hoed{J} \otimes J)$ of the right hand side to the other side,
we obtain
 \begin{eqnarray}
 \lefteqn{(J \otimes \hoed{J} \otimes J) W_{12} \Sigma_{12}
 (\hoed{J} \otimes \hoed{J} \otimes J) W_{12}^{\ast}W_{23}} \hspace{3cm} \nonumber \\
 &=& W_{23} (J \otimes \hoed{J} \otimes {J}) W_{12} \Sigma_{12}
 (\hoed{J} \otimes \hoed{J} \otimes J) W_{12}^{\ast}.  \hspace{2cm}  \label{eq-wj}
 \end{eqnarray}
By assumption, we have that
 $$\bolleke{W} = k \Sigma (\1 \otimes J \hoed{J} )
 \widetilde{W}^{\ast} W^{\ast} = k(J \otimes \hoed{J} )W \Sigma
 (\hoed{J} \otimes \hoed{J} ) W^{\ast}.$$
So equality~(\ref{eq-wj}) gives us that $W_{23}$ and $\bolleke{W}_{12}$ commute. Using the
fact that
 \bvgl
 J S^{\ast} J &=& \{(\omega \otimes \iota )((J \otimes J) \Sigma W^{\ast} \Sigma (J \otimes J))
 \mid \omega \in \B (\H )_{\ast} \} \\
 &=&  \{( \omega \otimes \iota ) (\bolleke{W} ) \mid \omega \in \B ( \H )_{\ast} \},
 \evgl
we conclude that $JMJ \subseteq M'$.

\item
We show that also the triple $(\hoed{W}, \hoed{J}, J)$ satisfies property (2'). Then it will
follow from the above that $\hoed{J} \hoed{M} \hoed{J} \subseteq \hoed{M}'$. We have
 \bvgl
 \lefteqn{\left((\1 \otimes J \hoed{J}\,) \Sigma (\hoed{W})\,{\bolleke{}} \; \hoed{W}
 (\hoed{W})\,\widetilde{} \right)^{\ast}} \hspace{3cm}\\
 &=& \Sigma (J \otimes J) \Sigma W^{\ast} \Sigma (J \otimes J) W (\hoed{J}
 \otimes \hoed{J}) \Sigma W^{\ast} \Sigma (\hoed{J} \otimes \hoed{J})(\1
 \otimes \hoed{J}J) \\
 &=& (\Sigma \bolleke{W} W \widetilde{W})(\1 \otimes \hoed{J}J)
 = k (\1 \otimes J \hoed{J}) (\1 \otimes \hoed{J}J) = k(\1 \otimes \1).
 \evgl
Hence
 $(\1 \otimes J \hoed{J}\,) \Sigma (\hoed{W})\,{\bolleke{}} \; \hoed{W}
 (\hoed{W})\,\widetilde{} = \overline{k} \1 \otimes \1.$

\item Finally we will show that $W$ is trim.  For $\omega \in \B(\H )_{\ast}$, we have that
 $$( \iota \otimes \omega) (\Sigma \widetilde{W})
 = (\iota \otimes \omega )(( \hoed{J} \otimes \hoed{J} ) W^{\ast}
 \Sigma (\hoed{J} \otimes \hoed{J}) );$$
the set of these elements is $\hoed{J} C^{\ast} \hoed{J}$. From the assumption it then follows
that
 \begin{equation} \label{eq-Cst}
 \hoed{J} C^{\ast} \hoed{J} = \{ ( \iota \otimes \omega )
 (\Sigma W^{\ast} \Sigma (J \otimes J)W(J \otimes \hoed{J})) \mid \omega
 \in \B(\H )_{\ast} \}.
 \end{equation}
Let $\xi, \eta \in \H$ and let $\omega$ denote the vector functional $\omega_{\xi, \eta}$.
Let $\theta \in \B(\H)_{\ast}$ be given by $\theta(x)=\omega(\hoed{J}x\hoed{J})^-$, so $\theta
= \omega_{\hoed{J}\xi,\hoed{J}\eta}$.  Like in the proof of lemma~\ref{lemma-set}, let $\{e_i
\mid i \in I\}$ be an orthonormal basis for $\H$.  Let $p_i$ denote the projection operator on
$\C e_i$.  We have
 \bvgl
 \lefteqn{(\iota \otimes \omega)(\Sigma W^{\ast}\Sigma (J \otimes J)
 W(J \otimes \hoed{J}))} \hspace{1cm} \\
 &=& J \; ((\iota \otimes \theta)
 (\Sigma W \Sigma (\1 \otimes \hoed{J}J)W)) \; J \\
 &=& J \; \left(\sum_{i \in I}(\iota \otimes \omega_{\hoed{J} \xi, \hoed{J} \eta})
 (\Sigma W \Sigma \; (\1 \otimes p_i) \; (\1 \otimes \hoed{J}J)W \right) \; J \\
 &=& J \; \left(\sum_{i \in I}(\iota \otimes \omega_{e_i, \hoed{J} \eta})(\Sigma W \Sigma)\;
 (\iota \otimes \omega_{\hoed{J}\xi,e_i})((\1 \otimes \hoed{J}J)W) \right) \; J \\
 &=& J \; \left(\sum_{i \in I}(\omega_{e_i, \hoed{J} \eta} \otimes \iota)(W)
 (\iota \otimes \omega_{\hoed{J} \xi,e_i}(\hoed{J}J \punt))(W) \right) \; J, \\
 \evgl
with convergence in (e.g.) the $\sigma$-weak topology.  This element belongs to the closure of
$J \hoed{M}MJ$.  Together with equality~(\ref{eq-Cst}), this gives us that
$\hoed{J}C^{\ast}\hoed{J}$ is contained in the closure of $J \hoed{M}MJ$. Therefore if $C$ is
$\sigma$-weakly dense in $\B(\H )$, then $W$ is trim.

\end{itemize}
\end{proof}

We can now show that $J$ and $\hoed{J}$ commute up to a scalar~:
\begin{proposition} \label{commute}
Let $(W,J,\hoed{J})$ be a quantum group frame.\\ There exists a complex number $\lambda$ with
modulus $1$ such that $J \hoed{J} = \lambda \hoed{J} J$.
\end{proposition}

\begin{proof}
Proposition \ref{prop-equivalent} tells us that there exists $k \in \C$ such that
 \begin{eqnarray}
 k (1 \otimes J \hoed{J} )
 &=& \Sigma \bolleke{W} W \widetilde{W} \nonumber \\
 &=& (J \otimes J) W^{\ast} (J \otimes J) \Sigma W \Sigma
 ( \hoed{J} \otimes \hoed{J} ) W^{\ast}
 ( \hoed{J} \otimes \hoed{J}) \Sigma \label{2}\\
 &=& (J \hoed{J} \otimes \1)W( \hoed{J} J \otimes \1)
 \Sigma W \Sigma (\1 \otimes \hoed{J} J)W(\1 \otimes J \hoed{J} ) \Sigma \label{3}.
 \end{eqnarray}
Taking adjoints of equality~(\ref{2}) we get
 \bvgl
 \lefteqn{\overline{k} (\1 \otimes \hoed{J} J)
 = \Sigma (\hoed{J} \otimes \hoed{J} )W( \hoed{J} \otimes
 \hoed{J}) \Sigma W^{\ast} \Sigma (J \otimes J)W(J \otimes J)}\\
 &=& \Sigma (J \otimes \hoed{J})(J \hoed{J} \otimes \1)W( \hoed{J} J \otimes \1)
 \Sigma W \Sigma (\1 \otimes \hoed{J} J)W(\1 \otimes J \hoed{J})
 \Sigma \Sigma (J \otimes \hoed{J}),
 \evgl
and using equality~(\ref{3}) we obtain
 $$\overline{k} (1 \otimes \hoed{J} J)= \Sigma (J \otimes \hoed{J} )
 k (\1 \otimes J \hoed{J} ) \Sigma (J \otimes \hoed{J} )
 = \overline{k} ( \hoed{J} J \hoed{J} J \otimes J \hoed{J}).$$
Hence $\hoed{J} J \hoed{J} J \otimes J \hoed{J} J \hoed{J} = \1$.  So there is $\lambda \in
\C$ such that $J \hoed{J}J \hoed{J} = \lambda \1$, and hence $J \hoed{J} = \lambda \hoed{J}
J$. Since both $J \hoed{J}$ and $\hoed{J}J$ are unitary operators, $\lambda$ must have modulus
$1$.
\end{proof}

Note that when $J$ is multiplied by a complex number $\gamma$ of modulus $1$, the properties
$J^2= \1$ and $J^{\ast} = J$ remain unchanged.  Also the anti-automorphism $R$ is not altered.
Hence if $(W, J, \hoed{J})$ is a quantum group frame, then so is $(W, \gamma J, \hoed{J})$.
Choosing $\gamma$ such that $\gamma^{-2} = \lambda$ gives us
 $$(\gamma J)\hoed{J} = \gamma \lambda \hoed{J}J
 = \overline{\gamma} \hoed{J}J = \hoed{J}(\gamma J).$$
This means that we have some choice when we define $J$, and we can determine $J$ in such a way
that $J$ and $\hoed{J}$ commute. However, often there are natural candidates for $J$ and
$\hoed{J}$, which do not necessarily commute
--- see the example of the $az+b$-group in section~\ref{sec-az+b}.  We will say something more
about this matter in remark~\ref{rem-sc}.

\lijntjetussen
As we mentioned in the introduction, the most mature theory of locally compact quantum groups
is the one introduced by Kustermans and Vaes in~\cite{kv2,kv3,kv4}.  In the following we argue
that the quantum groups which satisfy their definition, give rise to a quantum group frame.
Since we work in the von Neumann algebra setting, we will use the definition on the von
Neumann algebra level. Let us first recall this definition.

\bdefi {\textup(\cite{kv4} Definition 1.1)} \label{def-lcqg}
Consider a von Neumann algebra $M$ together with a unital normal $^*$-homomorphism $\Phi: M
\to M \otimes M$ such that $(\Phi \otimes \iota)\Phi = (\iota \otimes \Phi)\Phi$. Assume
moreover the existence of
\begin{enumerate}
\item
a normal semi-finite faithful (n.s.f.) weight $\varphi$ on $M$ that is left invariant, i.e.\
such that $\varphi((\omega \otimes \iota)\Phi(x))=\varphi(x) \omega(\1)$ for all $\omega \in
M^+_{\ast}$ and $x \in M^+$ such that $\varphi(x)<\infty$;
\item
a normal semi-finite faithful weight $\psi$ on $M$ that is right invariant, i.e.\ such that \\
$\psi ((\iota \otimes \omega)\Phi(x)) = \psi(x)\omega(\1)$ for all $\omega \in M^+_{\ast}$ and
$x \in M^+$ such that $\psi(x)<\infty$.
\end{enumerate}
Then the pair $(M,\Phi)$ is called a \emph{locally compact quantum group}.
\edefi

Let $(M, \Phi)$ be a locally compact quantum group.  Let $\varphi$ be a n.s.f.\ left invariant
weight on $(M,\Phi)$, and let $(\H, \iota, \Lambda)$ be a GNS-construction for $\varphi$.
Denote the modular conjugation of $\varphi$ with respect to this GNS-construction by $J$.
Theorem 2.1 of \cite{kv4} and the comments thereafter state that there exists a unique unitary
element $W \in \B(\H \otimes \H)$ such that
 $$W^{\ast}(\Lambda(x) \otimes \Lambda(y))
 = (\Lambda \otimes \Lambda)(\Phi(y)(x \otimes \1))$$
for all $x,y \in M$ such that $\varphi(x^{\ast}x)<\infty$ and $\varphi(y^{\ast}y)<\infty$.
This unitary satisfies the Pentagon equation. The $\sigma$-weak closure of the set $ S =
\{(\iota \otimes \omega)(W) \mid \omega \in \B(\H)_{\ast} \}$ is the von Neumann algebra $M$,
and for $x \in M$ one has that $\Phi(x) = W^{\ast}(\1 \otimes x)W$.  (Note that in~\cite{kv4}
the $\sigma$-strong $^{\ast}$-closure is considered; but since the norm closure of $S$ is a
$^*$-algebra, its $\sigma$-strong $^{\ast}$-closure coincides with the $\sigma$-weak closure.)

Define $\hoed{M}$ to be the $\sigma$-weak closure of the set $\{(\omega \otimes \iota)(W) \mid
\omega \in \B(\H)_{\ast} \}$. Then $\hoed{M}$ is a von Neumann algebra.  The correspondence
$\hoed{\Phi}(x)=\Sigma W (x \otimes \1) W^{\ast} \Sigma$ defines a unique normal
\mbox{$^*$-homomorphism} $\hoed{\Phi}:\hoed{M} \to \hoed{M} \otimes \hoed{M}$.  The pair
$(\hoed{M},\hoed{\Phi})$ is again a locally compact quantum group, called the \emph{dual} of
$(M,\Phi)$.  The dual left invariant weight $\hoed{\varphi}$ is constructed together with its
GNS-construction, having $\H$ as the underlying Hilbert space.  The modular conjugation of
$\hoed{\varphi}$ with respect to this GNS-construction is denoted as $\hoed{J}$.

\begin{proposition} \label{prop-frameqg}
The triple $(W,J,\hoed{J})$ is a quantum group frame on $\H$.
\end{proposition}

\begin{proof}
It follows from the theory of weights that $J$ and $\hoed{J}$ are isometric involutions of
$\H$ such that $M'=JMJ$ and $\hoed{M}'=\hoed{J} \hoed{M} \hoed{J}$.  The equality $(\hoed{J}
\otimes J)W(\hoed{J} \otimes J) = W^{\ast}$ is proved in corollary~2.2 of~\cite{kv4}.

From this equality we obtain that $\hoed{J}M\hoed{J}=M$ and $J\hoed{M}J=\hoed{M}$. One of the
comments after definition~1.6 of~\cite{kv4} says that $M \cap \hoed{M} = \C \1$.  Hence also
$\C \1 = J \hoed{J} \, M \, \hoed{J} J \; \cap \; J \hoed{J} \, \hoed{M} \, \hoed{J} J = M'
\cap \hoed{M}'$.  So $M \hoed{M} = \B(\H)$ and $W$ is trim.
\end{proof}

\bremark \label{rem-sc}
Denote the modular groups of $\varphi$ and $\psi$ by $\sigma$ and $\sigma'$ respectively. Then
proposition~6.8 of~\cite{kv2} gives us that there exists a strictly positive number $\nu$ such
that $\varphi \sigma'_t = \nu^t \varphi$ and $\psi \sigma_t = \nu^{-t} \psi$ for any $t \in
\R$.  The number $\nu$ is called the \emph{scaling constant}. Corollary 2.12 of~\cite{kv4} now
says that $\hoed{J}J=\nu^{i/4}J\hoed{J}$.  We remarked in proposition~\ref{commute} that $J$
could be rescaled in such a way that it commutes with $\hoed{J}$ and that $(W, J, \hoed{J})$
still is a quantum group frame.  But then $J$ is not anymore the anti-unitary obtained
naturally as the modular conjugation of $\varphi$.
\eremark

\section{Example~: The quantum $az+b$ group} \label{sec-az+b}

The quantum $az+b$ group was first introduced by Woronowicz in~\cite{wor-az+b}. The second
author continued the study of this example in~\cite{vd-az+b}.  He constructed right and left
Haar weights, and thus showed that it is a locally compact quantum group as in the
theory~\cite{kv2} of Kustermans and Vaes.  An important feature of this quantum group is that
the scaling constant (see remark~\ref{rem-sc}) is non-trivial.  This feature was foreseen by
the theory of Kustermans \& Vaes, and the quantum $ax+b$ group (see~\cite{wor-ax+b},
\cite{vd-az+b}) was the first example where it really occured.  The quantum $az+b$-group is
similar to the $ax+b$-group and has the same property, but is more simple. As we discussed in
remark~\ref{rem-sc}, in our theory this implies that the anti-unitary operators $J$ and
$\hoed{J}$ commute only up to a non-trivial scalar.

Since the quantum $az+b$-group is a locally compact quantum group, it follows from
proposition~\ref{prop-frameqg} that it gives rise to a quantum group frame.  The associated
multiplicative unitary and the anti-unitary $J$ are found in~\cite{vd-az+b}. In this section,
we will explicitely compute the anti-unitary operator $\hoed{J}$ and directly calculate the
commutation of $J$ and $\hoed{J}$.

\lijntjetussen
For the definition of the quantum $az+b$-group $(A,\Phi)$ and the notations used for its
ingredients, we refer to paragraph 3 of~\cite{vd-az+b}.  We will need to work with unbounded
operators acting on a Hilbert space; we refer to chapter 9 of \cite{stratila} for the basic
concepts and results.  We will not always be fully detailed when working with these unbounded
operators, but all the arguments of this section can be made rigorous.

Let $a$ and $b$ be as in definition~3.1 of~\cite{vd-az+b}, and let $a=u|a|$ and $b=v|b|$ be
their polar decompositions.  Consider the elements of the form
 $$\sum_{k,\ell = 0}^{2n-1} \left( \int f_{k,\ell}(|b|,t)|a|^{it} \,dt \right)
 v^k u^{\ell}$$
where the $f_{k,\ell}$ are continuous complex functions with compact support in $\R^+ \times
\R$, such that whenever $k\not= 0$ we have that $f_{k,\ell}(0,t) = 0$ for all $t$ and all
$\ell$. These elements form a non-degenerate $^*$-algebra $A_0$, which is norm dense in $A$
(\cite{vd-az+b} proposition 3.5 and definition 3.6).

The left Hilbert algebra associated to the right invariant weight on $(A,\Phi)$ is constructed
as follows (see proposition 4.2 of~\cite{vd-az+b}).   The underlying Hilbert space $\H$ is
$L^2(\R^+) \otimes L^2(\R) \otimes \C^{2n} \otimes \C^{2n}$, which will be identified with
$L^2(\R^+ \times \R, \C^{2n} \times \C^{2n})$. We fix an orthonormal basis $\{e_k \mid
k=0,1,\ldots,2n-1\}$ in $\C^{2n}$.  The index $k$ will be considered as an element of
$\Z_{2n}$; hence e.g.\ $e_{2n}$ stands for $e_0$. In what follows for $\xi \in \H$ we will use
the notation
 $$\xi = \sum_{k,\ell}\xi_{k,\ell} \otimes e_k \otimes e_{\ell}$$
with $\xi_{k,\ell} \in L^2(\R^+ \times \R)$.

When $x \in A_0$ is given by
 \begin{equation} \label{eq-x}
 x = \sum_{k,\ell}\left(\int f_{k,\ell}(|b|,t) \; |a|^{it} \, dt \right) v^k u^{\ell}
 \end{equation}
let $\eta(x) \in \H$ be given by
 $$(\eta(x))(r,t) = \sum_{k,\ell}e^{-\frac{\pi t}{n}} \; r^{\frac{1}{2}} \;
 f_{k,\ell}(r,t) \; e_k \otimes e_{\ell}$$
for $r \in \R^+$ and $t \in \R$.  The set of these $\eta(x)$ is made into a left Hilbert
algebra by letting
 $$\eta (x) \eta(y) = \eta(x y) \mbox{ and } \eta(x)^{\sharp} = \eta (x^{\ast})$$
whenever $x,y \in A_0$.

\begin{proposition}
The operator $J$ associated to the quantum \mbox{$az+b$} group is given by $J=J_{1} \otimes
J_{2}$, where $J_1$ acts on $L^2(\R^+ \times \R)$ as
 \bvgl
 (J_1 \xi)(r,t) &=& e^{-\frac{\pi t}{2n}} \; \overline{\xi}(e^{-\frac{\pi
 t}{n}}r,-t),
 \evgl
for $\xi \in K(\R^+ \times \R)$, $r \in \R^+$ and $t \in \R$, and $J_2$ acts on $\C^{2n}
\otimes \C^{2n}$ as
 \bvgl
 J_2(e_k \otimes e_{\ell}) &=& q^{k\ell} \; e_{-k} \otimes e_{-\ell}
 \evgl
for $k,\ell \in \{0,1,\ldots,2n-1\}$. The operator $\hoed{J}$ is given by $\hoed{J} =
\hoed{J}_1 \otimes \hoed{J}_2$ with
 $$\hoed{J}_1 = e^{\frac{- \pi i}{2n}}(e^{\frac{in}{2\pi}(\log a_0)^2} \otimes \1)
 \, e^{\frac{in}{\pi}( \log a_0 \otimes \log a_1)} \: C,$$
where $C$ denotes the complex conjugation on $L^2(\R^+ \times \R)$; and $J_2$ acts on
\mbox{$\C^{2n} \otimes \C^{2n}$} as
 $$\hoed{J}_2(e_k \otimes e_{\ell}) = (-1)^k \; q^{-\frac{1}{2} k^2} \; e_{-k} \otimes e_{k+\ell}$$
for $k,\ell \in \{0,1,\ldots,2n-1\}$.  Recall that $a_0$ and $a_1$ are the operators defined
as follows (cf. propositions 3.4 and 4.2 of~\cite{vd-az+b}).The operator $a_0$ on $L^2(\R^+)$
is defined by $(a_0^{it} f)(s) = e^{-\frac{\pi t}{2n}} f(e^{-\frac{\pi t}{n}} s)$ for $f \in
K(\R^+)$, and $a_1$ is the operator on $L^2(\R^)$ defined by $(a_1^{iq} \xi)(r) = \xi(r-q)$
for $\xi \in K(\R)$.
\end{proposition}

 \bproof
The operator $J$ is obtained in proposition 4.3 of~\cite{vd-az+b}.

The anti-unitary $\hoed{J}$ is to be obtained from the polar decomposition of the closure $G$
of the map $\eta(x) \mapsto \eta(\kappa(x)^{\ast})$, where $\kappa$ is the antipode of
$(A,\Phi)$ --- see proposition 2.8 and corollary 2.9 of~\cite{kv4}.  In proposition 3.12
of~\cite{vd-az+b}, the polar decomposition $\kappa = R \tau_{-\frac{i}{2}}$ of this antipode
is obtained. We will consider the map $\eta(x) \mapsto \eta(\kappa(x)^{\ast})$ as the
composition of the map $\eta(x) \mapsto \eta(\tau_{-\frac{i}{2}}(x))$ and the map $\eta(x)
\mapsto \eta(R(x)^{\ast})$.

We start with $\tau$.  Let $s \in \R$. Proposition~3.12 of~\cite{vd-az+b} gives us that
$\tau_s$ is defined by
 $$\begin{array}{lll}
 \tau_s(a) = a & \hspace{0.5cm} & \tau_s(b) = e^{\frac{2 \pi s}{n}} b.
 \end{array}$$
Hence
 $$\begin{array}{lll}
 \tau_s(|a|)= |a| & \hspace{0.5cm}  & \tau_s(|b|)  = e^{\frac{2 \pi s}{n}} |b|\\
 \tau_s(u) = u &  & \tau_s(v)  = v.
 \end{array}$$
Let $x \in A_0$ be given by~(\ref{eq-x}).  Then
 $$\tau_s(x)
 = \sum_{k,\ell}\left(\int f_{k,\ell}(e^{\frac{2 \pi s}{n}}|b|,t) \; |a|^{it} \, dt \right)
 v^k u^{\ell}$$
and
 \bvgl
 (\eta(\tau_s(x)))(r,t)
 &=& \sum_{k,\ell}  e^{-\frac{\pi t}{n}} \; r^{\frac{1}{2}} \;
     f_{k,\ell}(e^{\frac{2 \pi s}{n}}r,t) \; e_k \otimes e_{\ell}\\
 &=& \sum_{k,\ell} e^{-\frac{\pi s}{n}} \; e^{-\frac{\pi t}{n}} \;
     (e^{\frac{2 \pi s}{n}} r)^{\frac{1}{2}} \;
     f_{k,\ell}(e^{\frac{2 \pi s}{n}}r,t) \; e_k \otimes e_{\ell}.
 \evgl
Write $\xi=\eta(x)$ and $\widetilde{\xi} = \eta(\tau_s(x))$.  Then
 \bvgl
 \widetilde{\xi}_{k,\ell}(r,t)
 &=& e^{-\frac{\pi s}{n}} \; \xi_{k,\ell}(e^{\frac{2 \pi s}{n}}r,t)\\
 &=& e^{-\frac{2 \pi s}{n}} \; e^{\frac{\pi s}{n}} \;
 \xi_{k,\ell}(e^{\frac{2 \pi s}{n}}r,t)\\
 &=& e^{-\frac{2 \pi s}{n}} \; ((a_0^{-2is} \otimes \1) \xi_{k,\ell})(r,t),
 \evgl
where $a_0$ is given in proposition 3.4 of~\cite{vd-az+b}.  We conclude that
 $$ \widetilde{\xi} = e^{-\frac{2 \pi s}{n}}
 (a_0^{-2is} \otimes \1 \otimes \1 \otimes \1)\xi,$$
and so (the closure of) the map $\eta(x) \mapsto \eta(\tau_{-\frac{i}{2}}(x))$ is given by
 \begin{equation} \label{eq-tau}
 e^{\frac{\pi i}{n}}
 (a_0^{-1} \otimes \1 \otimes \1 \otimes \1).
 \end{equation}

Now we turn to the unitary antipode $R$.  We consider the map $x \mapsto R(x)^{\ast}$, which
is an anti-linear $^*$-automorphism.  Proposition~3.12 of~\cite{vd-az+b} gives us
 \begin{equation} \label{eq-RaRb}
 \begin{array}{lll}
 R(a)=a^{-1} &\hspace{0.5cm}& R(b)=-e^{\frac{\pi i}{n}}a^{-1}b.
 \end{array}
 \end{equation}
Hence we have that
 \bvgl
 R(a)^{\ast} &=& (a^{\ast})^{-1} = u \, |a|^{-1}\\
 R(b)^{\ast}
 &=& -e^{-\frac{\pi i}{n}} \, b^{\ast} \,(a^{\ast})^{-1}
 = - e^{-\frac{\pi i}{n}} \, v^{\ast} \, |b| \, u \, |a|^{-1}
 = - e^{-\frac{\pi i}{n}} \, v^{\ast} \, u \, |b| \, |a|^{-1}.
 \evgl
The operators $a$ and $b$ are defined in such a way that $|b|\,|a|^{-1} = e^{\frac{\pi i}{n}}
|a|^{-1} |b|$ (see definition 3.1 of~\cite{vd-az+b}). Therefore
 $$ R(b)^{\ast} = - v^{\ast} u \, |a|^{-1} \, |b|. $$
Like in~\cite{vd-az+b}, we write $q$ for the scalar $e^{\frac{\pi i}{n}}$.  For convenience
also the notations $q^{\frac{1}{2}}=e^{\frac{\pi i}{2n}}$ and $q^{-\frac{1}{2}}=e^{-\frac{\pi
i}{2n}}$ will be used.  We have
 $$ R(b)^{\ast} = - q^{-\frac{1}{2}} v^{\ast} u \, q^{\frac{1}{2}} |a|^{-1} \, |b|.$$
It follows from lemma~3.15 of~\cite{vd-az+b} that $- q^{-\frac{1}{2}} v^{\ast} u$ is unitary.
Lemma~3.14 of ~\cite{vd-az+b} gives us that (the closure of) $q^{\frac{1}{2}} |a|^{-1} \, |b|$
is a self-adjoint positive operator. We obtain
 $$\begin{array}{lll}
 R(u)^{\ast} = u  &\hspace{0.5cm}&
 R(v)^{\ast} =  - q^{-\frac{1}{2}} v^{\ast} u \\
 R(|a|)^{\ast} = |a|^{-1} &&
 R(|b|)^{\ast}= q^{\frac{1}{2}} |a|^{-1} \,|b|.
 \end{array}$$
Again let $x$ be as in~(\ref{eq-x}).  By the above, and since the map $x \mapsto R(x)^{\ast}$
is an anti-linear $^*$-automorphism, we have
 \begin{equation} \label{eq-R3}
 R(x)^{\ast}
 = \sum_{k,\ell} \left(\int \overline{f_{k,\ell}}
 (q^{\frac{1}{2}}|a|^{-1}\,|b|,t) |a|^{it} \, dt\right) (- q^{-\frac{1}{2}} v^{\ast} u)^k u^{\ell}.
 \end{equation}
By the definition of $a$ and $b$ we have $uv=qvu$; hence we have the commutation rule
 $$(v^{\ast}u)^k = v^{\ast}u v^{\ast} u \ldots v^{\ast}u
 = (q^{-1})^{\frac{1}{2}k(k-1)} (v^{\ast})^k u^k.$$
Therefore
 \begin{equation} \label{eq-R1}
 (- q^{-\frac{1}{2}} v^{\ast} u)^k = (-1)^k q^{-\frac{1}{2}k^2} (v^{\ast})^k u^k.
 \end{equation}
On the other hand, lemma~3.14 of \cite{vd-az+b} gives us that there is a unitary operator
$u_0$ such that $q^{\frac{1}{2}}|a|^{-1}\,|b| = u_0 |b| u_0^{\ast}$. This unitary $u_0$ is
given by $u_0= \exp (\frac{in}{2\pi} h^2)$ where $h=\log|a|$; hence it commutes with $|a|$. So
we have that
 \begin{equation} \label{eq-R2}
 \overline{f_{k,\ell}} (q^{\frac{1}{2}}|a|^{-1}\,|b|,t) \; |a|^{it}
 = u_0 \; \overline{f_{k,\ell}} (|b|,t) \; |a|^{it} \; u_0^{\ast}.
 \end{equation}
Substituting (\ref{eq-R1}) and (\ref{eq-R2}) in (\ref{eq-R3}), and taking into account that
$u_0$ commutes with $u$ and (by the definition of $a$ and $b$) also with $v$, we obtain that
 \bvgl
 R(x)^{\ast}
 &=& u_0 \left( \sum_{k,\ell} \left(\int \overline{f_{k,\ell}}
 (|b|,t) \, |a|^{it} \, dt \right) (-1)^k
 q^{-\frac{1}{2}k^2} v^{-k} u^{k+\ell} \right) u_0^{\ast}.
 \evgl
For $x \in A$, let $\pi(x)$ and $\pi'(x)$ be the GNS-operators given by
 $\pi(x) \eta(y) = \eta (xy)$ and $\pi'(x) \eta(y) = \eta (yx)$
if $y \in A_0$.  Then
 \begin{eqnarray*}
 \eta(R(x)^{\ast})
 &=& \pi(u_0) \: \eta \left( \left( \sum_{k,\ell} \left(\int \overline{f_{k,\ell}}
 (|b|,t) |a|^{it} \, dt \right)
 (-1)^k \; q^{-\frac{1}{2}k^2} \; v^{-k} \; u^{k+\ell} \right) u_0^{\ast} \right)\\
 &=& \pi(u_0) \: \pi'(u_0^{\ast}) \: \eta \left( \sum_{k,\ell} \left(\int \overline{f_{k,\ell}}
 (|b|,t) |a|^{it} \, dt \right)
 (-1)^k \; q^{-\frac{1}{2}k^2} \; v^{-k} \; u^{k+\ell} \right).
 \end{eqnarray*}
Hence the map $\eta(x) \mapsto R(x)^{\ast}$ is given by
 \begin{equation}
 \pi(u_0) \; \pi'(u_0^{\ast}) \; (C \otimes G_2) \label{eq-Rx}
 \end{equation}
where $C$ is the complex conjugation on $L^2(\R^+ \times \R)$, and $G_2$ is the anti-linear
operator on $L^2(\C^{2n} \times \C^{2n})$ given by
 $$G_2(e_k \otimes e_{\ell}) = (-1)^k \; q^{-\frac{1}{2} k^2} \; e_{-k} \otimes e_{k+\ell}.$$

We calculate $\pi(u_0) \pi'(u_0^{\ast})$. In proposition 4.2 of~\cite{vd-az+b}, the formula
 $$ \pi (a) = a_0 \otimes a_1 \otimes m \otimes s $$
is obtained. Here $m$ and $s$ are unitary, and $a_0$ and $a_1$ are positive operators.
Therefore, since $\pi$ is a $^*$-representation, it follows that
 $$ \pi (|a|) = a_0 \otimes a_1 \otimes \1 \otimes \1. $$
The situation for $\pi'$ is slightly more complicated, since it is an anti-representation but
not a $^*$-anti-representation.  Lemma 4.13 of~\cite{vd-az+b} gives us the formulas
 \bvgl
 \pi'(a) &=& q(\1 \otimes a_1 \otimes \1 \otimes s)\\
 \pi'(a^{\ast}) &=& q(\1 \otimes a_1 \otimes \1 \otimes s^{\ast}).\\
 \evgl
It follows that
 $$ \pi'(|a|) = q(\1 \otimes a_1 \otimes \1 \otimes \1).$$
Since there is no contribution in the last two legs of the tensor product, we forget about
these for a moment.  With $h = \log|a|$ we then have
 \bvgl
 \pi(h) &=& \log a_0 \otimes \1 + \1 \otimes \log a_1 \\
 \pi'(h) &=& \frac{\pi i}{n} \1 \otimes \1 + \1 \otimes \log a_1.
 \evgl
Therefore,
 \bvgl
 \pi(h^2) &=& (\log a_0)^2 \otimes \1 + \1 \otimes (\log a_1)^2
 + 2 \log a_0 \otimes \log a_1 \\
 \pi'(h^2) &=& - \frac{\pi^2}{n^2} \1 \otimes \1
 +\frac{2 \pi i}{n}(\1 \otimes \log a_1)
 +\1 \otimes (\log a_1)^2.
 \evgl
and with $u_0 = \exp (\frac{in}{2\pi} h^2)$ we obtain
 \bvgl
 \pi(u_0) &=& (e^{\frac{in}{2\pi}(\log a_0)^2}
 \otimes e^{\frac{in}{2\pi}(\log a_1)^2}) \,
 e^{\frac{in}{\pi}( \log a_0 \otimes \log a_1)} \\
 \pi'(u_0^{\ast}) &=& e^{-\frac{in}{2\pi} (- \frac{\pi^2}{n^2})} (\1 \otimes a_1)
 (\1 \otimes e^{-\frac{in}{2\pi}(\log a_1)^2}).
 \evgl
Hence
 \begin{equation} \label{eq-pi}
 \pi(u_0) \pi'(u_0^{\ast})
 = e^{\frac{\pi i}{2n}}(e^{\frac{in}{2\pi}(\log a_0)^2} \otimes a_1) \,
 e^{\frac{in}{\pi}( \log a_0 \otimes \log a_1)}.
 \end{equation}

We now bring together (\ref{eq-tau}), (\ref{eq-Rx}) and (\ref{eq-pi}).  We obtain
 $$G = e^{\frac{\pi i}{2n}}(e^{\frac{in}{2\pi}(\log a_0)^2} \otimes a_1)
 e^{\frac{in}{\pi}( \log a_0 \otimes \log a_1)} \: C \:
 e^{\frac{ \pi i}{n}} (a_0^{-1} \otimes \1)  \: \otimes G_2.$$
We see that the operators $\hoed{J}$ and $|G|$ in the polar decomposition $G=\hoed{J}\,|G|$
are given as follows~:
\begin{itemize}
\item $\hoed{J} = \hoed{J}_1 \otimes \hoed{J}_2$ with
 $$\hoed{J}_1 = e^{- \frac{\pi i}{2n}}(e^{\frac{in}{2\pi}(\log a_0)^2} \otimes \1)
 \, e^{\frac{in}{\pi}( \log a_0 \otimes \log a_1)} \: C,$$
and $\hoed{J}_2= G_2$;
\item $|G| = a_0^{-1} \otimes a_1 \otimes \1 \otimes \1$.
\end{itemize}

\eproof

\begin{proposition} \label{prop-commut}
$\hoed{J}J=e^{-\frac{\pi i}{n}} J \hoed{J}$.
\end{proposition}

\begin{proof}
We start with the anti-unitary operators $J_2$ and $\hoed{J}_2$.  For $k,\ell \in
\{0,1,\ldots,2n-1\}$ we have
 \bvgl
 J_2 \hoed{J}_2 (e_k \otimes e_{\ell})
 &=& J_2((-1)^k \; q^{-\frac{1}{2} k^2} \; e_{-k} \otimes e_{k+\ell})\\
 &=& (-1)^k \; q^{\frac{1}{2} k^2} \; q^{(-k)(k+\ell)} \; e_k \otimes e_{-k-\ell} \\
 &=& (-1)^k \; q^{- \frac{1}{2} k^2} \; q^{-k\ell} \; e_k \otimes e_{-k-\ell}.
 \evgl
On the other hand, we have
 \bvgl
 \hoed{J}_2 J_2 (e_k \otimes e_{\ell})
 &=& \hoed{J}_2(q^{k\ell} \; e_{-k} \otimes e_{-\ell})\\
 &=& q^{-k\ell} \; (-1)^k \; q^{-\frac{1}{2} (-k)^2} \; e_{k} \otimes e_{-k-\ell}.
 \evgl
So the operators $J_2$ and $\hoed{J}_2$ commute.

We now look at the operators $J_1$ and $\hoed{J}_1$.  Write $J_1=KC$ and
$\hoed{J}_1=\hoed{K}C$, where $C$ is the complex conjugation.  So $K$ is given by
 $$(K \xi)(r,t) = e^{-\frac{\pi t}{2n}} \xi(e^{-\frac{\pi t}{n}}r,-t)$$
for $\xi \in L^2(\R^+ \times \R)$ and $r \in \R^+$, $t \in \R$.  Note that $KC=CK$.  Write
$K=(\1 \otimes K')K''$, where $K'$ is the unitary operator on $L^2(\R)$ given by $(K' \xi)(t)
= \xi(-t)$, and $K''$ is the unitary operator on $L^2(\R^+ \times \R)$ given by
 $$(K'' \xi)(r,t)=
 e^{\frac{\pi t}{2n}} \xi(e^{\frac{\pi t}{n}}r,t) = ((a_0^{-it}\otimes \1)\xi)(r,t).$$
The operator $\hoed{K}$ is given by
 $$\hoed{K} = e^{- \frac{\pi i}{2n}}(e^{\frac{in}{2\pi}(\log a_0)^2} \otimes \1)
 e^{\frac{in}{\pi}( \log a_0 \otimes \log a_1)}.$$
We calculate $C \hoed{K} C$.  Since $C a_0 C = a_0^{-1}$ and $C a_1 C = a_1^{-1}$, we have
 $$C \hoed{K}C = e^{\frac{\pi i}{2n}}(e^{- \frac{in}{2\pi}(\log a_0)^2} \otimes \1)
 e^{- \frac{in}{\pi}( \log a_0 \otimes \log a_1)}.$$

We know that $J_1$ and $\hoed{J}_1$ commute up to a scalar factor.  Since these operators have
their first leg in the commutative \cst-subalgebra of $\B(L^2(\R^+))$ generated by
$\log(a_0)$, we can find this scalar factor by evaluating the first leg of $J_1 \hoed{J}_1$
and $\hoed{J}_1 J_1$ in any point of the spectrum of $a_0$. Formally, this corresponds to
replacing $\log a_0$ by a scalar $s$.  Then $\hoed{K}$ corresponds to $e^{- \frac{\pi i}{2n}}
\; e^{\frac{ins^2}{2\pi}} \; a_1^{\frac{ins}{\pi}}$, and $C \hoed{K} C$ corresponds to
$e^{\frac{\pi i}{2n}} \, e^{- \frac{ins^2}{2\pi}} \, a_1^{ \frac{ins}{\pi}}$.  Let $k$ denote
the self-adjoint operator on $L^2(\R)$ given by $k \xi (t) = t\xi(t)$ for $\xi \in K(\R)$ and
$t \in \R$.  Then $K''$ corresponds to $e^{-isk}$ and $K=(\1 \otimes K')K''$ corresponds to
$K' e^{-isk}$.  Hence $J_1 \hoed{J}_1 = K(C \hoed{K} C)$ corresponds to
 \begin{equation} \label{eq-20}
 K' \, e^{-isk} \, e^{\frac{\pi i}{2n}} \, e^{- \frac{ins^2}{2\pi}} \,
 a_1^{ \frac{ins}{\pi}}.
 \end{equation}
The operator $\hoed{J}_1 J_1 = \hoed{K} C K C = \hoed{K} K$ corresponds to
 \begin{eqnarray}
 e^{- \frac{\pi i}{2n}} \; e^{\frac{ins^2}{2\pi}} \;
    a_1^{\frac{ins}{\pi}} \, K' \, e^{-isk}
 &=& e^{- \frac{\pi i}{2n}} \; e^{\frac{ins^2}{2\pi}} \;
    K' a_1^{- \frac{ins}{\pi}} \, e^{-isk} \nonumber \\
 &=& K' e^{- \frac{\pi i}{2n}} \; e^{\frac{ins^2}{2\pi}} \;
    a_1^{- \frac{ins}{\pi}} \, e^{-isk}, \label{eq-vr}
 \end{eqnarray}
where we used $a_1^{\frac{ins}{\pi}} K' = K' a_1^{-\frac{ins}{\pi}} $. Now for $p \in \R$ we
have
 $$(a_1^{ip} \; e^{-isk} \, f)(t) = (e^{-isk} \; f)(t-p) = e^{-is(t-p)} \;
 f(t-p),$$
while
 $$(e^{-isk} \, a_1^{ip} \; f)(t) = (e^{-ist} \; f)(t-p);$$
so $a_1^{ip} \; e^{-isk} = e^{isp} \; e^{-isk} \, a_1^{ip}$.  Using this in~(\ref{eq-vr}) we
obtain that $\hoed{J}_1J_1$ corresponds to
 \begin{equation} \label{eq-21}
 K' e^{- \frac{\pi i}{2n}} \; e^{\frac{ins^2}{2\pi}} \; e^{is(-\frac{ns}{\pi})}
    \, e^{-isk} a_1^{- \frac{ins}{\pi}}.
 \end{equation}
Comparing (\ref{eq-20}) and (\ref{eq-21}) we obtain that $\hoed{J}_1J_1=e^{-\frac{\pi i}{n}}
J_1 \hoed{J}_1$.  Since $J_2$ and $\hoed{J}_2$ commute, this gives us the result.
 \eproof

\bremark
In remark 3.3.3, we mentioned that in the setting of the left regular representation
$\hoed{J}J = \nu^{i/4} J \hoed{J}$, where $\nu$ is the scaling constant. Recall how the
translation between left and right approach can be made.  The $(W,J,\hoed{J})$ associated to
the right Haar measure we considered in this example, is in fact the ''left'' quantum group
frame corresponding to the dual of the quantum $az+b$ group with the opposite
comultiplication. Passing to the dual with the opposite comultiplication does not alter the
scaling constant. The first author showed in~\cite{vd-az+b} that this scaling constant equals
$e^{-\frac{4 \pi}{n}}$. So proposition~\ref{prop-commut} corresponds to the result predicted
by the theory.
\eremark

\section{Crossed product with a locally compact group}

In classical group theory, the construction of the semi-direct product is a powerful tool to
construct classes of examples of groups.  In this section we generalize this procedure, to
obtain the crossed product of a quantum group frame with a locally compact group.

Let $G,H$ be locally compact groups, and let $\alpha$ be a continuous action of $G$ on $H$. We
look at this situation on the group von Neumann algebra level. The group action gives rise to
an action of $G$ on $\VN(H)$, which we will denote by $\widetilde{\alpha}$. It is the
homomorphism of $G$ into the group of $^*$-automorphisms of $\VN(H)$ given by
 $$\widetilde{\alpha}_g (\lambda_{h}) = \lambda_{\alpha_{g}(h)}$$
if $g \in G$ and $h \in H$.  One can then construct a unitary representation $u$ of $G$ on
$L^2(H)$, which implements this action, i.e.
 $$\widetilde{\alpha}_{g} (m)= u_{g} m u_{g}^{\ast}$$
for any $g \in G$ and $m \in \M(\H)$.

We denote the quantum group frame corresponding to the group von Neumann algebra $\VN(H)$ by
$(W^H, J^H, \hoed{J^H})$. Note that this is the situation dual to the one discussed in
example~\ref{ex-group2}; we use superscripts here.  The fact that $\alpha$ acts as
automorphisms of $H$, translates into the property that for any $g \in G$, the operator $u_g
\otimes u_g$ commutes with $W^H$, and $u_g$ commutes with $J^H$ and $\hoed{J^H}$.

These data can be generalized to the quantum group frame setting as follows.

\bdefi \label{defaction}
Let $\H$ be a Hilbert space and $(W,J, \hoed{J} )$ a quantum group frame on $\H$. Let $G$ be a
locally compact group and
 $$u:G \rightarrow \B (\H):p \mapsto u_{p}$$
a (strongly continuous) unitary representation of $G$ on $\H$. We say that $u$ is an
\emph{action} of $G$ on $(W,J, \hoed{J})$ if for any $p \in G$, we have that
\begin{enumerate}
\item $W( u_{p} \otimes u_{p} )=( u_{p} \otimes u_{p})W,$
\item $J u_{p} = u_{p} J$ and $\hoed{J} u_{p}  = u_{p} \hoed{J}$.
\end{enumerate}
\edefi

The group $G$ will come into play in its form of group von Neumann algebra.  We will denote
the associated quantum group frame by $(W^0, J^0, \hoed{J^0})$. For further reference, we
write down the formulas once again~:
 \begin{eqnarray}
 &&W^0 \xi (p,q) = \xi (q p,q)\\
 &&J^0 \eta (p) = \Delta (p)^{-\frac{1}{2}} \; \overline{\eta (p^{-1})}\\
 &&\hoed{J^0} \eta (p) = \overline{\eta (p)}
 \end{eqnarray}
whenever $\xi \in L^2 (G \times G)$, $\eta \in L^2(G)$ and $p,q \in G$.

We lift the action $u$ to the level of the operator algebras. Recall (see
e.g.~\cite{boekfons}) that the space $K(G, \H)$ of $\H$-valued continuous functions on $G$
with compact support can be identified with a dense subspace of $\H \otimes L^2(G)$.  For $\xi
\in K(G,\H)$, let $U \xi: G \rightarrow \H$ be given by
 $$(U \xi)(p) = u_{p^{-1}} \xi (p)$$
if $p \in G$.  Note that $U\xi$ is again an element of $K(G,\H)$, which has the same norm as
$\xi$.  Hence $U$ extends to a unitary operator on $\H \otimes L^2(G)$.

\bdefi \label{def-crossprod}
Define the operators $W^1 \in \B ( \H \otimes L^2 (G) \otimes \H \otimes L^2 (G))$ and $J^1,
\hoed{J^1}$ on $\H \otimes L^2 (G)$ by
\begin{itemize}
\item $W^1 = W_{24}^0 U_{14} W_{13} = U_{14}W^0_{24}W_{13}$,
\item $J^1 = (J \otimes J^0 )U = U^{\ast}(J \otimes J^0)$,
\item $\hoed{J^1} = \hoed{J} \otimes \hoed{J^0}$.
\end{itemize}
\lijntjeomhoog \edefi

\bremark \label{rem-comm}
The equality of operators
 \begin{equation} \label{comm1}
 W_{23}^0 U_{13} = U_{13} W_{23}^0.
 \end{equation}
is easily checked on the continuous $\H$-valued functions on $G$ with compact support. Since
these functions form a dense subspace of the Hilbert space $\H \otimes L^2(G)$, we conclude
that $W_{23}^0$ and $U_{13}$ indeed commute. Similarly, one can verify the relation
 \begin{equation}
 \label{comm5} (J \otimes J^0)U=U^{\ast}(J \otimes J^0).
 \end{equation}
\eremark

Note that in the group example the formulas of the definition above give us the quantum group
frame of (the group von Neumann algebra of) the semi-direct product of $H$ and $G$ by the
action $\alpha$.

The formulas of definition~\ref{def-crossprod} also correspond to the formulas of the crossed
product of a Kac system with an automorphism group, which Baaj \& Skandalis construct in
proposition 8.16 of~\cite{bs}.  Those formulas are based on the right regular representation;
the translation to the ''left regular'' setting is made by noting that when $W$ is the left
regular representation of a quantum group, the unitary $\Sigma W^{\ast} \Sigma$ corresponds to
the right regular representation of the quantum group with the opposite comultiplication.

\lijntjetussen
In general, we have the following result.  The necessary background on the crossed product of
a von Neumann algebra with a locally compact group can be found in e.g.~\cite{boekfons}.

\begin{proposition} \label{prop-crossedproduct}
Let $(W, J, \hoed{J})$ be a quantum group frame on a Hilbert space $\H$.  Let $G$ be a locally
compact group, and $u$ an action of $G$ on $(W, J, \hoed{J})$. Let $W^1$, $J^1$ and
$\hoed{J^1}$ be defined as in definition~\ref{def-crossprod}. Then $(W^1, J^1, \hoed{J^1})$ is
a quantum group frame on $\H \otimes L^2 (G)$. The corresponding von Neumann algebra $M^1$ is
the crossed product $M \cross_{\alpha} G$, where $\alpha$ is the action of $G$ on $M$ given by
$\alpha_p (m) = u_p m u_p^{\ast}$ if $m \in M$ and $p \in G$. The dual von Neumann algebra
$\hoed{M^1}$ is the (ordinary) von Neumann algebra tensor product $\hoed{M} \otimes
L^{\infty}(G)$.
\end{proposition}

\begin{proof}
\begin{enumerate}
\item
In the same way as those in remark~\ref{rem-comm}, the following commutation relations are
proved~:
\begin{eqnarray}
 \label{comm2}
 W_{23}^0 U_{12} &=& U_{12} U_{13} W_{23}^0;\\
 \label{comm3}
 U_{13} U_{23} &=& U_{23} U_{13};\\
 \label{comm4}
 U_{13} U_{23} W_{12} &=& W_{12} U_{13} U_{23};\\
 \label{comm6}
 (\hoed{J} \otimes J^0)U &=& U^{\ast}(\hoed{J} \otimes J^0).
 \end{eqnarray}

\item
Clearly, $W^1$ is unitary. When we refer to $(\H \otimes L^2(G))^{\otimes 3}$ as a three-fold
tensor product, we will use brackets for the leg numbering. When we consider it as a six-fold
tensor product we use the leg numbering notation as before.  By the equalities~(\ref{comm1}),
(\ref{comm4}) and again (\ref{comm1}), we have that
 \bvgl
 W_{(12)}^1 W_{(13)}^1 W_{(23)}^1
 &=& U_{14} W_{24}^0 W_{13} U_{16} W_{26}^0 W_{15} U_{36} W_{46}^0 W_{35}\\
 &=& U_{14} W_{24}^0 W_{26}^0 W_{46}^0 W_{13} U_{16} U_{36} W_{15} W_{35}\\
 &=& U_{14} W_{24}^0 W_{26}^0 W_{46}^0 U_{16} U_{36} W_{13} W_{15} W_{35}\\
 &=& U_{14} U_{16} W_{24}^0 W_{26}^0 W_{46}^0 U_{36} W_{13} W_{15} W_{35}.
 \evgl
Using the multiplicativity of $W^0$ and of $W$, and the commutation relation~(\ref{comm2}),
this gives us
 $$
 W_{(12)}^1 W_{(13)}^1 W_{(23)}^1
 = W_{46}^0 U_{14} W_{24}^0 U_{36} W_{35} W_{13}\\
 = W_{(23)}^1 W_{(12)}^1.
 $$
So $W^1$ is multiplicative.

\item
Now we will show that $M^1$ is the crossed product $M \cross_{\alpha} G$.  Let $\omega \in \B
( \H )_{\ast}$ and $\omega^0 \in \B (L^2 (G))_{\ast}$. Then, if we denote $( \iota \otimes
\omega)(W)$ by $m$, we have that
 $$(\iota \otimes (\omega \otimes \omega^0))(W^1)
 = ( \iota \otimes \iota \otimes \omega \otimes \omega^0 ) ( U_{14} W_{24}^0 W_{13} )
 = (\iota \otimes \iota \otimes \omega^0)(U_{13} W_{23}^0)\;(m \otimes \1).$$
Denoting $(\iota \otimes \omega^0)(W^0)$ by $m^0$, we obtain
 $$
 (\iota \otimes \iota \otimes \omega^0)(U_{13} W_{23}^0)
 = U^{\ast}((\iota \otimes \iota \otimes \omega^0)(U_{12} U_{13} W_{23}^0))\\
 = U^{\ast}(\1 \otimes m^0)U,
 $$
where we used equality~(\ref{comm2}) in the last step.  Observe that for $p \in G$, we have
that $U^{\ast}(\1 \otimes \lambda_p) U = u_p \otimes \lambda_p$.  So $M^1$ is the von Neumann
algebra generated by the sets
 $$ \{m \otimes \1 \mid m \in M \} \: \mbox{ and } \:
 \{u_p \otimes \lambda_p \mid p \in G \}.$$
Hence $M^1$ is spatially isomorphic to the crossed product $M \cross_{\alpha} G$~: \newline in
fact $M^1 = \mbox{$U^{\ast} (M \cross_\alpha G)U$}$ (see e.g.~proposition 2.12
of~\cite{boekfons}).

\item
In this part we show that $\hoed{M^1}=\hoed{M} \otimes L^{\infty}(G)$.  We will for a moment
work on the \mbox{\cst-algebraic} level.  When $A$ is a \cst-algebra, we denote its multiplier
algebra (see e.g.~\cite{wor-pseudo}) by $M(A)$. Recall that $W^0 \in M(C^{\ast}_r(G) \otimes
C_0(G))$, where $C^{\ast}_r(G)$ denotes the reduced group \cst-algebra of $G$ and $C_0(G)$
denotes the \cst-algebra of continuous complex functions on $G$ vanishing at infinity.
Therefore the operator $W^0$ can be considered as a strictly continuous bounded
$M(C^{\ast}_r(G))$-valued function on $G$. We will denote the \cst-algebra of compact
operators on a Hilbert space $\H$ as $\K(\H)$. The unitary operator $U$ belongs to the
multiplier algebra $M(\K(\H) \otimes C_0(G))$, and hence can be considered as a strictly
continuous bounded function from $G$ to $M(\K(\H))$.

First note that using the fact that $W^0 \in \M(G) \otimes L^{\infty}(G)$, that $U \in
M(\K(\H) \otimes C_0(G)) \subseteq \B(\H) \otimes L^{\infty}(G)$ and that $W \in M \otimes
\hoed{M}$, we obtain
 $$W^1 = W^0_{24} U_{14} W_{13} \in \B(\H) \otimes \M(G) \otimes \hoed{M} \otimes
 L^{\infty}(G),$$
and hence $\hoed{M^1} \subseteq \hoed{M} \otimes L^{\infty}(G)$. We need to show that this
last inclusion is an equality.

On the \cst-algebraic level, from the observations above, it follows that
 $$W^1 \in M(\K(\H) \otimes C^{\ast}_r(G) \otimes \K(\H) \otimes C_0(G)).$$
Hence $\hoed{S^1} \subseteq M(\K(\H) \otimes C_0(G))$, and we can consider elements of
$\hoed{S^1}$ as $\B(\H)$-valued functions on $G$.  Under this identification, we have for
$\theta \in \B (\H)_{\ast}$ and $\theta^0 \in \B (L^2 (G))_{\ast}$ that
 \begin{eqnarray}
 \lefteqn{((\theta \otimes \theta^0 \otimes \iota \otimes \iota)(W^1))(p)
 =((\1 \otimes (\theta^0 \otimes \iota)(W^0))(\theta \otimes \iota \otimes
 \iota)(U_{13} W_{12}))(p)} \hspace{1cm} \nonumber \\
 && = \: ((\theta^0 \otimes \iota)(W^0))(p)\:
 (\theta (u_{p^{-1}} \punt) \otimes \iota)(W) \hspace{2.8cm}  \label{vergelijking1}
 \end{eqnarray}
for any $p \in G$.

Let $\xi, \eta \in \H$ and $\xi^0, \eta^0 \in K(G)$ with support in a compact $K \subseteq G$.
Consider the vector functionals $\omega = \omega_{\xi, \eta}$ on $\H$ and $ \omega^0 =
\omega_{\xi^0, \eta^0}$ on $L^2 (G)$. We will prove that
 $$(\omega \otimes \iota)(W) \otimes (\omega^0 \otimes
 \iota)(W^0) \in \hoed{M^1},$$
and as linear combinations of such elements are dense in $\hoed{M} \otimes L^{\infty} (G)$,
this will provide the inclusion $\hoed{M} \otimes L^{\infty}(G) \subseteq \hoed{M^1}$.

Let $\eps > 0$ be arbitrary.  For each $q \in K$ we fix an open neighbourhood $E_q$ of $q$
such that
 $$\| u_{p^{-1}} \eta - u_{q^{-1}} \eta \| < \frac{\eps}{\| \xi^0 \| \:
 \| \eta^0 \| \: \| \xi \|}$$
for any $p \in E_q$. This is possible because the map $G \rightarrow \H: p \mapsto u_{p^{-1}}
\eta$ is continuous. In this way we obtain an open cover $\{ E_q \mid q \in K \}$ of $K$. Let
$ \mbox{$\mathcal{E}$} = \{ E_{p_1}, E_{p_2}, \ldots, E_{p_n} \}$ be a finite subcover.

Let $A(G)$ denote the Fourier algebra of the group $G$ (see \cite{eymard} chapitre 3). In
lemma~\ref{lemma-partitie} below we will show that it is possible to take a partition of unity
$\{ h_1, h_2, \ldots, h_n \}$ of K subordinate to $\mathcal{E}$, in such a way that $h_i \in
K(G) \cap A(G)$ for every $i \in \{1,2,\ldots,n\}$. Recall (\cite{eymard} th\'{e}or\`{e}me
3.10 and following comments) that $A(G)$ is isomorphic to the von Neumann algebra predual of
$\M(G)$. We will identify these two algebras, and for $h \in A(G)$ denote the corresponding
$\sigma$-weakly continuous linear functional on $\M(G)$ again by $h$. The duality is given by
 $$h \left(\int f(p) \lambda_p \,dp \right) = \int h(p) f(p) \,dp$$
whenever $h \in A(G)$ and $f \in L^1(G)$.  For $i \in \{1,2,\ldots,n\}$, define $\omega_i$ to
be $h_i \punt \omega^{0}$ as an element of $\M(G)_{\ast}$. So $\omega_i$ is the linear
functional on $\VN(G)$ given by
 $$\omega_i(x) = (h_i \punt \omega^0)(x)
 = (h_i \otimes \omega^0)(\Phi^0(x))$$
if $x \in \VN(G)$, where $\Phi^0$ is the comultiplication on $\VN(G)$ induced by $W^0$. As a
strictly continuous $M(C^{\ast}_r)$-valued function on $G$, we have that $W^0$ is given by
$W^0(p) = \lambda_p$.  Hence
 \begin{eqnarray}
 ((\omega_i \otimes \iota)(W^0))(p) &=& (h_i \punt \omega^0)(\lambda_p) \\
 &=& h_i(p)\omega^0(\lambda_p) \nonumber \\
 &=& h_i(p)((\omega^0 \otimes \iota)(W^0))(p). \label{vergelijking2}
 \end{eqnarray}
As both $h_i$ and $\omega^0$ are $\sigma$-weakly continuous linear functionals on $\VN(G)$, so
will be their product $\omega_i$. Therefore it can be extended to a $\sigma$-weakly continuous
linear functional on $\B(L^2(G))$, which we again will denote by $\omega_i$. From
equality~(\ref{vergelijking1}), it follows for $p \in G$ that
 $$((\omega(u_{p_i}\punt) \otimes \omega_i \otimes
 \iota \otimes \iota)(W^1))(p) = ((\omega_i \otimes
 \iota)(W^0))(p)\;(\omega(u_{p_i}u_{p^{-1}}\punt) \otimes \iota)(W).$$
Using~(\ref{vergelijking2}), we obtain that
 \bvgl
 \lefteqn{ \left\| \left( \sum_{i=1}^n
 (\omega (u_{p_i} \punt) \otimes \omega_i \otimes \iota \otimes \iota)(W^1)\right) (p)
 - (\omega \otimes \iota)(W)\,((\omega^0 \otimes \iota) (W^0))(p) \right\|  }\\
 &=& |((\omega^0 \otimes \iota)(W^0)) (p) | \,
 \left\| \sum_{i=1}^n h_i (p) ((\omega (u_{p_i p^{-1}} \punt) \otimes \iota )(W)
 - (\omega \otimes \iota )(W)) \right\|\\
 &\leq& \| \xi^0 \| \, \| \eta^0 \| \, \sum_{i=1}^{n} h_i(p) \,
 \| \omega (u_{p_i} u_{p^{-1}} \punt) - \omega \|\\
 &\leq& \| \xi^0 \| \, \| \eta^0 \| \,
 \sum_{i=1}^{n} h_i (p) \, \| \xi \| \, \| u_{p} u_{p_i^{-1}}\eta - \eta \| \\
 &\leq& \|\xi^0 \| \, \| \eta^0 \| \, \| \xi \| \,
 \sum_{i=1}^n h_i(p) \| u_{p} \|  \, \| u_{p_i^{-1}} \eta - u_{p^{-1}} \eta \|.
 \evgl
Now $h_i(p)=0$ unless $p \in E_i$, in which case
 $$\| u_{p^{-1}} \eta - u_{p_i^{-1}} \eta \|
 < \frac{\eps}{\| \xi^0 \| \, \| \eta^0 \| \, \| \xi\|}.$$
So the above norm is less than
 $$\| \xi^0 \| \, \| \eta^0 \| \, \| \xi \| \,
 \sum_{i=1}^{n} h_i (p) \frac{\eps}{\| \xi^0 \| \, \| \eta^0\| \, \| \xi \|}
 = \eps.$$
We have obtained that $(\omega \otimes \iota)(W) \otimes (\omega_0 \otimes \iota)(W^0)$,
considered as a $\B(\H)$-valued functions on $G$,  can be uniformly approximated by elements
of $\hoed{S^1}$, where $\B(\H)$ is considered with its norm topology. This uniform topology is
stronger than the $\sigma$-weak topology on $\hoed{S^1}$. Therefore we can conclude that
$(\omega \otimes \iota)(W) \otimes (\omega_0 \otimes \iota)(W_0)$ belongs to the $\sigma$-weak
closure of $\hoed{S^1}$, and hence to $\hoed{M^1}$.

\item
From parts (3)\ and (4)\ of this proof, it follows that the $\sigma$-weak closures of $S^1$
and $\hoed{S^1}$ are self-adjoint.  Therefore, these $\sigma$-weak closures are $M^1$ and
$\hoed{M^1}$ respectively. We have obtained that
 \bvgl
 M^1 &=& \{m \otimes \1, u_p \otimes \lambda_p \mid m \in M, p \in G \}'', \\
 \hoed{M^1} &=& \hoed{M} \otimes L^{\infty}(G).
 \evgl
It follows that $M^1 \hoed{M^1}$ contains $M \hoed{M} \otimes \1$.  Since $W$ is trim, $M^1
\hoed{M^1}=\B(\H)$ and it follows that $M^1 \hoed{M^1}$ also contains $\1 \otimes
\VN(G)\,L^{\infty}(G)$.  Hence from the fact that $W$ and $W^0$ are trim, we obtain that also
$W^1$ is trim.

\item
It is clear that $J^1=(J^1)^{\ast}$, $(J^1)^2 = \1$, $\hoed{J^1}= (\hoed{J^1})^{\ast}$ and
$(\hoed{J^1})^2 = \1$. We will show that $(\hoed{J^1} \otimes J^1)W^1 (\hoed{J^1} \otimes
J^1)=(W^1)^{\ast}$. By equality~(\ref{comm6}), we have
 \bvgl
 (\hoed{J^1} \otimes J^1)W^1(\hoed{J^1} \otimes J^1)
 &=& U_{34}^{\ast} U_{14}^{\ast}
 (\hoed{J} \otimes \hoed{J^0} \otimes J \otimes J^0) W^0_{24} W_{13}
 ( \hoed{J} \otimes \hoed{J^0} \otimes J \otimes J^0) U_{34} \\
 &=& U_{34}^{\ast} U_{14}^{\ast} (W_{24}^0)^{\ast} W_{13}^{\ast} U_{34}.
 \evgl
Using the commutation relations~(\ref{comm4}) and (\ref{comm3}), we obtain that this is
 $$W_{13}^{\ast} U_{34}^{\ast} U_{14}^{\ast} U_{34} (W_{24}^0)^{\ast} =
 W_{13}^{\ast} U_{14}^{\ast} (W_{24}^0)^{\ast} = (W^1)^{\ast}.$$

\item
The von Neumann algebra $$J^1 M^1 J^1 = U^{\ast} (J \otimes J^0)M^1(J \otimes J^0)U$$ is
generated by the elements
\begin{itemize}
\item $U^{\ast}(J \otimes J^0)(m \otimes \1)(J \otimes J^0)U = U^{\ast}
(JmJ \otimes \1)U$ where $m \in M$, and
\item $(J \otimes J^0)U U^{\ast} (\1 \otimes \lambda_p)U U^{\ast}
(J \otimes J^0)=\1 \otimes J^0 \lambda_p J^0$ where $p \in G$.
\end{itemize}
By proposition 3.11 of~\cite{boekfons} and equality~(\ref{comm5}), we have that
 \bvgl
 \1 \otimes J^0 \lambda_p J^0
 &=& (J \otimes J^0)U(u_p \otimes \lambda_p)U^{\ast}(J \otimes J^0)\\
 &=& U^{\ast} (J u_p J \otimes J^0 \lambda_p J^0) U
 = U^{\ast} (\1 \otimes \rho_p)U,
 \evgl
where $\rho_p$ denotes the operator of right translation given by
 $$(\rho_p \xi )(q)= \Delta (p)^{\frac{1}{2}} \xi (qp)$$
for $\xi \in L^2(G)$ and $q \in G$.

Using the fact that $JMJ \subseteq M'$, we obtain that $J^1 M^1 J^1$ is contained in the von
Neumann algebra generated by
 $$\{ U^{\ast} (m' \otimes \1)U \mbox{, } U^{\ast} (u_p \otimes \rho_p)U
 \mid m' \in M', p \in G \}$$
which (by e.g.\ theorem 3.12 of~\cite{boekfons}) is exactly $U^{\ast} (M \cross_{\alpha}
G)'U=(M^1)'$. With $\hoed{M^1} = \hoed{M} \otimes \hoed{M^0}$ we immediately have that
$\hoed{J^1} \hoed{M^1} \hoed{J^1} \subseteq \hoed{M}' \otimes (\hoed{M^0})'=(\hoed{M^1})'$.

\end{enumerate}
\end{proof}

We still need to show that a partition of unity on a locally compact group $G$ can be chosen
in $K(G) \cap A(G)$.

\begin{lemma} \label{lemma-partitie}
Let $G$ be a locally compact group.
 Let $K$ be a compact subset of $G$ and $\mathcal{E}= \{ E_1, E_2, \ldots,
E_n \}$ a finite open cover of $K$.  Then there are positive functions $\{ h_1, h_2, \ldots,
h_n \}$ in $A(G) \cap K(G)$ such that
\begin{enumerate}
\item  $0 \leq \sum_{i=1}^n h_i \leq 1$,
\item  $\sum_{i=1}^n h_i(p) = 1$ for any $p \in K$,
\item  For any $i \in \{1,2,\ldots n\}$, the function $h_i$ has its support in $E_i$.

\end{enumerate}
\end{lemma}

\begin{proof}
We denote the support of a function $f$ by $\mbox{supp}(f)$. Let $D$ be an open neighbourhood
of the unit $e$ of $G$ with compact closure $\overline{D}$, such that $K \overline{D}
\subseteq E_1 \cup E_2 \cup \ldots \cup E_n$.  Let $\{ f_1, f_2, \ldots, f_n \}$ be a
partition of unity of $K \overline{D}$ subordinate to $\mathcal{E}$.  Fix an open
neighbourhood $E$ of $e$ such that $E \subseteq D$ and $(\mbox{supp}(f_i))E^{-1} \subseteq
E_i$ for all $i$.  Let $g \in K(G)$ such that $g$ has its support in $E$, $g\geq 0$ and $\int
g(p)\,dp = 1$.

For $i \in \{1,2,\ldots n\}$ and $p \in G$, let
 $$ h_i(p) = \int g(p^{-1}q) f_i(q) \,dq.$$
For $f,g \in K(G)$, we denote by $\widetilde{f}$ the function in $K(G)$ given by
 $$\widetilde{f}(p) = \overline{f(p^{-1})},$$
and by $f \star g$ the function given by
 $$(f \star g)(p) = \int f(q) g(q^{-1}p) \, dq.$$
Then we can write the above equality as
 $$ h_i(p) = \int f_i(q) \widetilde{g}(q^{-1}p) \,dq = (f_i \star \widetilde{g})(p).$$
Since for $p \in G$, we have that $(f_i \star \tilde{g})(p)=\langle \lambda_{p^{-1}} f_i \mid
g \rangle$, the functions $f_i \star \tilde{g}$ will belong to the predual $A(G)$ of $\VN(G)$.
Hence the functions $h_i$ are in $A(G)$.  Since both $g$ and $f_i$ are positive elements of
$K(G)$, also $h_i \in K(G)$ and is positive.

\begin{enumerate}
\item
We have that \bvgl 0 &\leq& \sum_{i=1}^n h_i(p) = \int g(p^{-1}q) \sum_{i=1}^n f_i(q) \,dq \\
&\leq& \int g(p^{-1}q) \,dq = \int g(q) \,dq = 1. \evgl

\item
Let $p \in K$.  We have that $$\sum_{i=1}^n h_i(p) = \int g(p^{-1}q) \sum_{i=1}^n f_i(q)
\,dq.$$  Now $g(p^{-1}q) \not= 0$ is only possible if $p^{-1}q \in E$, and hence $q \in pE
\subseteq p \overline{D}$.  Since $\{f_1, f_2, \ldots, f_n \}$ is a partition of unity of $K
\overline{D}$, it follows that $\sum_{i=1}^n f_i(q)=1$. So we have that $g(p^{-1}q)
\sum_{i=1}^n f_i(q) = g(p^{-1}q)$ for any $q \in G$. We conclude that $\sum_{i=1}^n h_i(p) =
1$.

\item
Let $i \in \{1,2,\ldots n\}$ and $p \in K$ such that $h_i(p) \not= 0$.  Then there must be a
$q \in G$ such that $g(p^{-1}q) \not= 0$ and $f_i(q) \not= 0$. So $p^{-1}q$ must be in $E$,
and $q$ must be an element of $\mbox{supp}(f_i)$. Then $p \in qE^{-1} \subseteq (\mbox{supp}
(f_i))E^{-1} \subseteq E_i$.  We conclude that $h_i$ has its support in $E_i$.

\end{enumerate}
\end{proof}

\bremark  \label{rem-def} Let us take a closer look at the von Neumann algebras $M^1$ and
$\hoed{M^1}$. The \emph{algebra structure} of $M^1$ is a crossed product. Both $M$ and $M^0
=\VN(G)$ are embedded in $M^1$, but not in the trivial way as in $M \otimes M^0$. We could say
that the algebra structure of $M \otimes M^0$ has been twisted by the given action. On the
other hand, the \emph{coalgebra structure} is trivially given by the embeddings $m \mapsto m
\otimes \1$ of $M$ and $\lambda_p \mapsto U^{\ast} (\1 \otimes \lambda_p)U$ of $\VN(G)$  in
$M^1$. For $m \in M$ we indeed have
 \bvgl
 \Phi^1 (m \otimes \1) &=& (W^1)^{\ast} (\1 \otimes \1 \otimes m \otimes \1)W^1 \\
 &=& W_{13}^{\ast} (W^0)_{24}^{\ast} U_{14}^{\ast}
 (\1 \otimes \1 \otimes m \otimes \1)U_{14} W_{24}^0 W_{13}\\
 &=& \Phi(m)_{13},
 \evgl
and for $p \in G$ we have
 \bvgl
 \Phi^1 (U^{\ast} (\1 \otimes \lambda_p)U)
 &=& (W^1)^{\ast} U_{34}^{\ast} (\1 \otimes \1 \otimes \1 \otimes
 \lambda_p)U_{34} W^1 \\
 &=& W_{13}^{\ast} (W^0)_{24}^{\ast}U_{14}^{\ast} U_{34}^{\ast}
 (\1 \otimes \1 \otimes \1 \otimes \lambda_p) U_{34} U_{14}W_{24}^0 W_{13}\\
 &=& U_{14}^{\ast} U_{34}^{\ast} W_{13}^{\ast} (W^0)_{24}^{\ast}
 (\1 \otimes \1 \otimes \1 \otimes \lambda_p) W_{24}^0 W_{13} U_{34} U_{14} \\
 &=& U_{14}^{\ast} U_{34}^{\ast} (\1 \otimes \lambda_p \otimes \1 \otimes \lambda_p)
 U_{34} U_{14}\\
 &=& (U^{\ast} \otimes U^{\ast}) \Phi (\lambda_p)_{24} (U \otimes U),
 \evgl
where in the third step we used the commutation relations (\ref{comm4}) and (\ref{comm1}). For
the dual $\hoed{M^1}$, we have that the \emph{algebra structure} trivially is the tensor
product $\hoed{M} \otimes \hoed{M^0}$. But now the \emph{coalgebra structure} is deformed~:
For $\hoed{m} \in \hoed{M}$ and $f \in \hoed{M^0}= L^{\infty} (G)$, we have that
 \bvgl
 \hoed{\Phi^1} (\hoed{m} \otimes f)
 &=& (\hoed{W^1})^{\ast} (\1 \otimes \1 \otimes \hoed{m} \otimes f) \hoed{W^1} \\
 &=& U_{32} (\hoed{W^0})_{24}^{\ast} \hoed{W}_{13}^{\ast}
 (\1 \otimes \1 \otimes \hoed{m} \otimes f) \hoed{W}_{13} (\hoed{W^0})_{24} U_{32}^{\ast} \\
 &=& (\iota \otimes \tau \otimes \iota)(\hoed{\Phi} \otimes \hoed{\Phi^0})
 (\hoed{m} \otimes f),
 \evgl
where $\tau$ denotes the twisted flip
 $$\tau: \B(\H \otimes L^2(G)) \rightarrow \B(L^2(G) \otimes \H):
 x \otimes x^0 \mapsto \sigma(U(x \otimes x^0)U^{\ast}).$$
\end{remark}

\end{document}